\documentclass[MSNbibl,number,citesort,seceqn,dvips]{arxbj}
\usepackage{accents}

\aid{0}
\volume{19}
\issue{3}
\pubyear{2013}
\firstpage{846}
\lastpage{885}
\doi{10.3150/12-BEJ468} 

\makeatletter

\newcommand{\RMe}{\mathrm{e}}

\newcommand{\TV}{\mathrm{TV}}

\renewcommand{\underbar}{\underaccent{\bar}}

\renewcommand{\mid}{|}

\newcommand{\rrvert}{\vert}
\newcommand{\llvert}{\vert}

\newcommand{\RMo}{\mathrm{o}}

\newcommand{\mrmd}{\,\mathrm{d}}

\newproclaim{definition}{Definition}[section]
\newtheorem{theorem}[definition]{Theorem}
\newtheorem{lemma}[definition]{Lemma}
\newtheorem{proposition}[definition]{Proposition}
\newtheorem{corollary}[definition]{Corollary}

\newremark{remark}[definition]{Remark}

\newcommand{\pinf}{p_{\mathrm{inf}}}
\newcommand{\Prob}{\operatorname{Pr}}

\makeatother

\begin{document}
\begin{frontmatter}

\title{Divergence rates of Markov order estimators and their
application to statistical estimation of stationary ergodic processes}
\runtitle{Statistical estimation of stationary ergodic processes}

\begin{aug}
\author{\fnms{Zsolt} \snm{Talata}\corref{}\ead[label=e1]{talata@math.ku.edu}}
\runauthor{Zs. Talata} 
\address{Department of Mathematics, University of Kansas, 405 Snow Hall,
1460 Jayhawk Boulevard, Lawrence, KS 66045-7523, USA. \printead{e1}}
\end{aug}

\received{\smonth{4} \syear{2010}}
\revised{\smonth{6} \syear{2011}}

%
\begin{abstract}
Stationary ergodic processes with finite alphabets are estimated by
finite memory processes from a sample, an $n$-length realization of the
process, where the memory depth of the estimator process is also
estimated from the sample using penalized maximum likelihood (PML).
Under some assumptions on the continuity rate and the assumption of
non-nullness, a rate of convergence in $\bar{d}$-distance is obtained,
with explicit constants. The result requires an analysis of the
divergence of PML Markov order estimators for not necessarily finite
memory processes. This divergence problem is investigated in more
generality for three information criteria: the Bayesian information
criterion with generalized penalty term yielding the PML, and the
normalized maximum likelihood and the Krichevsky--Trofimov code
lengths. Lower and upper bounds on the estimated order are obtained.
The notion of consistent Markov order estimation is generalized for
infinite memory processes using the concept of oracle order estimates,
and generalized consistency of the PML Markov order estimator is
presented.
\end{abstract}

%

\end{frontmatter}

\section{Introduction}

This paper is concerned with the problem of estimating stationary
ergodic processes with finite alphabet from a sample, an observed
length $n$ realization of the process, with the $\bar{d}$-distance
being considered between the process and the estimated one. The $\bar
{d}$-distance was introduced by Ornstein~\cite{O1} and became one of
the most widely used metrics over stationary processes. Two stationary
processes are close in $\bar{d}$-distance if there is a joint
distribution whose marginals are the distributions of the processes
such that the marginal processes are close with high probability (see
Section~\ref{secappl} for the formal definition). The class of ergodic
processes is $\bar{d}$-closed and entropy is $\bar{d}$-continuous,
which properties do not hold for the weak topology~\cite{ShB}.

Ornstein and Weiss~\cite{OW} proved that for stationary processes
isomorphic to i.i.d. processes, the empirical distribution of the
$k(n)$-length blocks is a strongly consistent estimator of the
$k(n)$-length parts of the process in $\bar{d}$-distance if and only if
$k(n)\le(\log n)/h$, where $h$ denotes the entropy of the
process.

Csisz\'ar and Talata~\cite{CsT3} estimated the $n$-length part of a
stationary ergodic process $X$ by a Markov process of order
$k_n$. The transition probabilities of this Markov
estimator process are the empirical conditional probabilities, and the
order $k_n\to+\infty$ does not depend on the sample. They obtained a
rate of convergence of the Markov estimator to the process $X$ in
$\bar{d}$-distance, which consists of two terms. The first one is the
bias due to the error of the approximation of the process by a Markov
chain. The second term is the variation due to the error of the
estimation of the parameters of the Markov chain from a
sample.

In this paper, the order $k_n$ of the Markov estimator process is
estimated from the sample. For the order estimation, penalized maximum
likelihood (PML) with general penalty term is used. The resulted Markov
estimator process finds a tradeoff between the bias and the variation
as it uses shorter memory for faster memory decays of the process $X$.
If the process $X$ is a Markov chain, the PML order estimation recovers
its order asymptotically with a wide range of penalty terms.

Not only an asymptotic rate of convergence result is obtained but also
an explicit bound on the probability that the $\bar{d}$-distance of the
above Markov estimator from the process $X$ is greater than
$\varepsilon
$. It is assumed that the process $X$ is non-null, that is, the
conditional probabilities of the symbols given the pasts are separated
from zero, and that the continuity rate of the process $X$ is summable
and the restricted continuity rate is uniformly convergent. These
conditions are usually assumed in this area~\cite{Br,DGG,FG,Marton}.
The summability of the continuity rate implies that the process is
isomorphic to an i.i.d. process~\cite{Bb}.

The above result on statistical estimation of stationary ergodic
processes requires a non-asymptotic analysis of the Markov order
estimation for not necessarily finite memory processes. In this paper,
this problem is also investigated in more generality: under milder
conditions than it would be needed for the above bound and not only for
the PML method.

A popular approach to the Markov order estimation is the minimum
description length (MDL) principle~\cite{RissanenB,BRY}. This method
evaluates an information criterion for each candidate order based on
the sample and the estimator takes the order for which the value is
minimal. The normalized maximum likelihood (NML)~\cite{St} and the
Krichevsky--Trofimov (KT)~\cite{KT} code lengths are natural
information criteria because the former minimizes the worst case
maximum redundancy for the model class of $k$-order Markov chains,
while the latter does so, up to an additive constant, with the average
redundancy. The Bayesian information criterion (BIC)~\cite{Schw} can be
regarded as an approximation of the NML and KT code lengths. The PML is
a generalization of BIC; special settings of the penalty term yield the
BIC and other well-known information criteria, such as the Akaike
information criterion (AIC)~\cite{Akaike}. There are other methods for
Markov order estimation, see~\cite{Narayan} and references there, and
the problem can also be formulated in the setting of hypothesis testing
\cite{Ryabko}.

If a process is a Markov chain, the NML and KT Markov order estimators
are strongly consistent if the candidate orders have an upper bound
$\RMo(\log n)$~\cite{Cs}. Without such a bound, they fail to be consistent
\cite{CsSh}. The BIC Markov order estimator is strongly consistent
without any bound on the candidate orders~\cite{CsSh}. If a process has
infinite memory, the Markov order estimators are expected to tend to
infinity as $n\to+\infty$. The concept of context trees of arbitrary
stationary ergodic processes is a model more complex than Markov
chains. Recent results~\cite{CsT} in that area imply that this
expectation holds true for the BIC and KT Markov order estimators but
they provide no information about the asymptotics of the
divergence.

In this paper, the divergence of the PML, NML and KT Markov order
estimators for not necessarily finite memory processes is
investigated.
Not only asymptotic rates of divergence are obtained but also explicit
bounds on the probability that the estimators are greater and less,
respectively, than some order. Instead of the usual assumption of
non-nullness, it is assumed only that the conditional probabilities of
one of the symbols given the pasts are separated from zero. This
property is called weakly non-nullness and is ``noticeably weaker''
than non-nullness~\cite{CFF}.

First, the process is assumed to be weakly non-null and $\alpha
$-summable. The $\alpha$-summability~\cite{N1,N2,GL,L} is a condition
weaker than the summability of the continuity rate. Under these
conditions, a bound on the probability that the estimators are greater
than some order is obtained, that yields an $\mathcal{O}(\log n)$ upper
bound on the estimated order eventually almost surely as $n\to+\infty$.

Then, a bound on the probability that the estimators are less than some
order is obtained assuming that the process is weakly non-null and the
decay of its continuity rates is in some exponential range. This bound
implies that the estimators satisfying the conditions attain a $c\log
n$ divergence rate eventually almost surely as $n\to+\infty$, where the
coefficient $c$ depends on the range of the continuity rates. The class
of processes with exponentially decaying continuity rate is considered
in various problems~\cite{DGG,GGG}. Fast divergence rate of the
estimators are expected only for a certain range of continuity rates.
Clearly, the estimators do not have a fast divergence rate if the
memory decay of the process is too fast. On the other hand, too slow
memory decay is also not favored to a fast divergence rate because then
the empirical probabilities do not necessarily converge to the true
probabilities.

To provide additional insight into the asymptotics of Markov order
estimators, the notion of consistent Markov order estimation is
generalized for infinite memory processes. A Markov order estimator is
compared to its oracle version, which is calculated based on the true
distribution of the process instead of the empirical distribution. The
oracle concept is used in various problems, see, for example,
\cite{Barron,Birge,Ooracle,GaoGijbels}. If the decay of the continuity
rate of the process is faster than exponential, the ratio of the PML
Markov order estimator with sufficiently large penalty term to its
oracle version is shown to converge to $1$ in probability.

The structure of the paper is the following. In Section \ref
{secnotation}, notation and definitions are introduced for stationary
ergodic processes with finite alphabets. In Section~\ref{secic}, the
PML, NML and KT information criteria are introduced. Section \ref
{secmain} contains the results on divergence of the
information-criterion based Markov order estimators. In Section \ref
{secappl}, the problem of estimating stationary ergodic process in
$\bar{d}$-distance is formulated and our results are presented. The
results require bounds on empirical entropies, which are stated in
Section~\ref{secmain} and are proved in Section~\ref{secent}.
Section~\ref{secproof} contains the proof of the divergence results,
and Section~\ref{secproofappl} the proof of the process estimation results.

\section{Finite and infinite memory processes}\label{secnotation}

Let $X= \{ X_i, -\infty<i<+\infty\}$ be a stationary ergodic
stochastic process with finite alphabet $A$. We write $X_i^j =
X_i,\ldots,X_j$ and $x_i^j = x_i,\ldots,x_j \in A^{j-i+1}$ for $j\ge
i$. If $j<i$, $x_i^j$~is the empty string. For two strings $x_1^i\in
A^i$ and $y_1^j\in A^j$, $x_1^i y_1^j$ denotes their concatenation
$x_1,\ldots,x_i,y_1,\ldots,y_j \in A^{i+j}$. Write
\[
P\bigl(x_i^j\bigr) = \Prob\bigl( X_i^j
= x_i^j \bigr)
\]
and, if $P(x_{-m}^{-1})>0$,
\[
P\bigl(a | x_{-m}^{-1} \bigr) = \Prob\bigl(
X_0 = a \mid X_{-m}^{-1} =
x_{-m}^{-1} \bigr).
\]
For $m=0$, $P(a | x_{-m}^{-1} ) = P ( a)$.

The process $X$ is called \textit{weakly non-null} if
\[
\alpha_0 = \sum_{a\in A}
\inf_{x_{-\infty}^{-1} \in A^{\infty}} P\bigl(a | x_{-\infty}^{-1} \bigr
) >0.
\]
Letting
\[
\alpha_k = \min_{y_{-k}^{-1} \in A^k} \sum_{a\in A}
\inf_{x_{-\infty
}^{-1} \in A^{\infty}: x_{-k}^{-1} = y_{-k}^{-1}} P\bigl(a | x_{-\infty
}^{-1} \bigr),\qquad k=1,2,
\ldots,
\]
we say that the process $X$ is \textit{$\alpha$-summable} if
\[
\alpha= \sum_{k=0}^{+\infty} (1-
\alpha_k) <+\infty.
\]

The \textit{continuity rates} of the process $X$ are
\[
\bar{\gamma}(k) = \sup_{x_{-\infty}^{-1} \in A^{\infty}} \sum_{a\in A}
\bigl\llvert P\bigl(a | x_{-k}^{-1} \bigr) - P\bigl(a |
x_{-\infty}^{-1} \bigr) \bigr\rrvert
\]
and
\[
\underbar{\gamma}(k) = \inf_{x_{-\infty}^{-1} \in A^{\infty}} \sum
_{a\in A}
\bigl\llvert P\bigl(a | x_{-k}^{-1} \bigr) - P\bigl(a |
x_{-\infty}^{-1} \bigr) \bigr\rrvert.
\]
Obviously, $\underbar{\gamma}(k) \le\bar{\gamma}(k)$. If $\sum
_{k=1}^{\infty} \bar{\gamma}(k)<+\infty$, then the process $X$ is said
to have \textit{summable continuity rate}.
%
%
\begin{remark}\label{remgammaeq}
Since for any $x_{-k}^{-1} \in A^k$ and $z_{-m}^{-k-1}\in A^{m-k}$,
$m\ge k$,
\[
\inf_{x_{-\infty}^{-k-1}} P\bigl(a | x_{-\infty}^{-1} \bigr) \le P
\bigl(a | z_{-m}^{-k-1} x_{-k}^{-1} \bigr)
\le\sup_{x_{-\infty}^{-k-1}} P\bigl(a | x_{-\infty}^{-1} \bigr),
\]
the above definition of continuity rate is equivalent to
\[
\bar{\gamma}(k) = \sup_{ i>k } \max_{x_{-i}^{-1} \in A^i} \sum
_{a\in A} \bigl\llvert P\bigl(a | x_{-k}^{-1}
\bigr) - P\bigl(a | x_{-i}^{-1} \bigr) \bigr\rrvert.\
\]
\end{remark}
%
%
\begin{remark}\label{remgammaalpha}
The process is $\alpha$-summable if it has summable continuity rate because
\begin{eqnarray*}
1-\alpha_k
&\le&1 - \max_{y_{-k}^{-1} \in A^k} \sum_{a\in A} P
\bigl(a | y_{-k}^{-1} \bigr)\\
&&{} + \max_{y_{-k}^{-1} \in A^k} \sum
_{a\in A} \sup_{x_{-\infty}^{-1}
\in
A^{\infty}: x_{-k}^{-1} = y_{-k}^{-1}} \bigl( P\bigl(a |
y_{-k}^{-1} \bigr) - P\bigl(a | x_{-\infty}^{-1}
\bigr) \bigr)
\\
&\le& |A| \bar{\gamma}(k).
\end{eqnarray*}
\end{remark}

The $k$-order \textit{entropy} of the process $X$ is
\[
H_k = - \sum_{a_1^k\in A^k} P
\bigl(a_1^k\bigr) \log P\bigl(a_1^k
\bigr),\qquad k\ge1,
\]
and the $k$-order \textit{conditional entropy} is
\[
h_k = - \sum_{a_1^{k+1}\in A^{k+1}} P
\bigl(a_1^{k+1}\bigr) \log P\bigl(a_{k+1} |
a_1^k\bigr),\qquad k\ge0.
\]
Logarithms are to the base $2$. It is well known for stationary
processes~\cite{Cover,CSbook} that the conditional entropy $h_k$ is a
non-negative decreasing function of $k$, therefore its limit exists as
$k\to+\infty$. The \textit{entropy rate} of the process is
\[
\bar{H} = \lim_{k\to+\infty} h_k = \lim_{k\to+\infty} \frac1k
H_k.
\]
Note that $h_k - \bar{H} \ge0$ for any $k\ge0$.

The process $X$ is a \textit{Markov chain} of order $k$ if for each $n>k$
and $x_1^n\in A^n$
%
%
\begin{equation}
\label{eqMCdef} P \bigl(x_1^n\bigr) = P
\bigl(x_1^k\bigr) \prod_{i=k+1}^n
P\bigl(x_i | x_{i-k}^{i-1}\bigr),
\end{equation}
where $P (x_1^k)$ is called initial distribution and $ \{
P(a|a_1^k), a\in A, a_1^k\in A^k \}$ is called transition
probability matrix. The case $k=0$ corresponds to i.i.d. processes. The
process $X$ is of \textit{infinite memory} if it is not a Markov chain
for any order $k<+\infty$. For infinite memory processes, $h_k - \bar
{H}>0$ for any $k\ge0$.

In this paper, we consider statistical estimates based on a sample
$X_1^n$, an $n$-length part of the process. Let $N_n (a_1^k)$ denote
the number of occurrences of the string $a_1^k$ in the sample $X_1^n$
\[
N_n \bigl(a_1^k\bigr) = \bigl\llvert\bigl
\{ i: X_{i+1}^{i+k} = a_1^k, 0\le i
\le n-k \bigr\}\bigr\rrvert.
\]
For $k\ge1$, the empirical probability of the string $a_1^k$ is
\[
\hat{P} \bigl(a_1^k\bigr) = \frac{N_n (a_1^k)}{n-k+1}
\]
and the empirical conditional probability of $a_{k+1}\in A$ given
$a_1^k$ is
\[
\hat{P} \bigl(a_{k+1}| a_1^k\bigr) =
\frac{ N_n(a_1^{k+1}) }{ N_{n-1}(a_1^k) }.
\]
For $k=0$, $\hat{P} (a_{k+1}| a_1^k) = \hat{P} (a_{k+1})$.
The $k$-order \textit{empirical entropy} is
\[
\hat{H}_{k}\bigl(X_1^n\bigr) = - \sum
_{a_1^k\in A^k} \hat{P} \bigl(a_1^k
\bigr) \log\hat{P} \bigl(a_1^k\bigr),\qquad 1\le k\le n,
\]
and the $k$-order \textit{empirical conditional entropy} is
\[
\hat{h}_{k}\bigl(X_1^n\bigr) = - \sum
_{a_1^{k+1}\in A^{k+1}} \hat{P} \bigl(a_1^{k+1}
\bigr) \log\hat{P} \bigl(a_{k+1}| a_1^k\bigr),\qquad
0\le k\le n-1.
\]

The likelihood of the sample $X_1^n$ with respect to a $k$-order Markov
chain model of the process $X$ with some transition probability matrix
$ \{ Q(a_{k+1}|a_1^k), a_{k+1}\in A, a_1^k\in A^k \}$, by~(\ref{eqMCdef}), is
\[
P'\bigl( X_1^n \bigr) = P'
\bigl(X_1^k\bigr) \prod_{a_1^{k+1}\in A^{k+1}}
Q\bigl(a_{k+1}| a_1^k\bigr)^{N_n(a_1^{k+1})}.
\]
For $0\le k< n$, the \textit{maximum likelihood} is the maximum in
$Q(a_{k+1}|a_1^k)$ of the second factor above, which equals
\[
\mathrm{ML}_{k}\bigl(X_1^n\bigr) = \prod
_{a_1^{k+1}\in A^{k+1}} \hat{P} \bigl(a_{k+1}|
a_1^k\bigr)^{N_n(a_1^{k+1})}. 
\]
Note that $\log\mathrm{ML}_{k}(X_1^n) = -(n-k) \hat{h}_{k}(X_1^n)$.

\section{Information criteria}\label{secic}

An information criterion assigns a score to each hypothetical model
(here, Markov chain order) based on a sample, and the estimator will be
that model whose score is minimal.
%
%
\begin{definition}\label{defIC}
For an information criterion
\[
\mathrm{IC}_{X_1^n} ( \cdot)\dvtx\mathbb{N}\to\mathbb{R}^{+},
\]
the Markov order estimator is
\[
\hat{k}_{\mathrm{IC}} \bigl(X_1^n\bigr) = \arg
\min_{0\le k < n} \mathrm{IC}_{X_1^n} (k).
\]
\end{definition}
%
%
\begin{remark}
Here, the number of candidate Markov chain orders based on a sample is
finite, therefore the minimum is attained. If the minimizer is not
unique, the smallest one will be taken as $\arg\min$.
\end{remark}

We consider three, the most frequently used information criteria,
namely, the Bayesian information criterion and its generalization, the
family of penalized maximum likelihood (PML) \mbox{\cite{Schw,CsSh}}, the
normalized maximum likelihood (NML) code length~\cite{St}, and the
Krichevsky--Trofimov (KT) code length~\cite{KT}.
%
%
\begin{definition}\label{defPML}
Given a penalty function $\operatorname{pen}(n)$, a non-decreasing
function of the
sample size $n$, for a candidate order $0\le k<n$ the PML criterion is
\begin{eqnarray*}
\mathrm{PML}_{X_1^n} ( k ) &=& - \log\mathrm{ML}_{k}
\bigl(X_1^n\bigr) + \bigl(|A|-1\bigr) |A|^k
\operatorname{pen}(n)
\\
&=& (n-k) \hat{h}_{k}\bigl(X_1^n\bigr) +
\bigl(|A|-1\bigr) |A|^k \operatorname{pen}(n).
\end{eqnarray*}
\end{definition}

The $k$-order Markov chain model of the process $X$ is described by the
conditional probabilities $ \{ Q(a_{k+1}|a_1^k), a_{k+1}\in A,
a_1^k\in A^k \}$, and $(|A|-1) |A|^k$ of these are free parameters.

The second term of the PML criterion, which is proportional to the
number of free parameters of the $k$-order Markov chain model, is
increasing in $k$. The first term, for a given sample, is known to be
decreasing in $k$. Hence, minimizing the criterion yields a tradeoff
between the goodness of fit of the sample to the model and the
complexity of the model.
%
%
\begin{remark}
If $\operatorname{pen}(n) = \frac12\log n$, the PML criterion is
called \textit{Bayesian information criterion} (BIC), and if
$\operatorname{pen}(n)
= 1$, \textit{Akaike information criterion} (AIC).
\end{remark}

The minimum description length (MDL) principle minimizes the length of
a code of the sample tailored to the model class. Strictly speaking,
the information criterion would have an additive term, the length of a
code of the structure parameter. This additional term, the length of a
code of $k$, is omitted since it does not affect the results.
%
%
\begin{definition} \label{defNML}
For a candidate order $0\le k<n$, the NML criterion is
\[
\mathrm{NML}_{X_1^n} ( k ) = - \log P_{\mathrm{NML}, k}
\bigl(X_1^n\bigr),
\]
where
\[
P_{\mathrm{NML}, k} \bigl(X_1^n\bigr) =
\frac{ \mathrm{ML}_{k}(X_1^n) }{
\Sigma(n,k) } \qquad\mbox{with } \Sigma(n,k) = \sum
_{x_1^n\in A^n} \mathrm{ML}_{k} \bigl(x_1^n
\bigr)
\]
is the $k$-order NML-probability of $X_1^n$.
\end{definition}
%
%
\begin{remark}
Writing
\[
\mathrm{NML}_{X_1^n} ( k ) = - \log\mathrm{ML}_{k}
\bigl(X_1^n\bigr) + \log\Sigma_{n,k},
\]
the NML criterion can be regarded as a PML criterion in a broader sense.
\end{remark}
%
%
\begin{definition} \label{defKT}
For a candidate order $0\le k<n$, the KT criterion is
\[
\mathrm{KT}_{X_1^n} ( k ) = - \log P_{\mathrm{KT}, k}
\bigl(X_1^n\bigr),
\]
where
\[
P_{\mathrm{KT}, k} \bigl(X_1^n\bigr) =
\frac{ 1 }{ |A|^k } \!\!\mathop{\prod_{a_1^k\in A^k:}}_{N_{n-1}(a_1^k)\ge1}\!\!
\frac{
\prod_{a_{k+1}: N_n(a_1^{k+1}) \ge1} [
( N_n(a_1^{k+1}) - 1/2 ) ( N_n(a_1^{k+1}) -
3/2 )
\cdots( 1/2 ) ]
}{
( N_{n-1}(a_1^k) - 1 + {|A|}/{2} ) (
N_{n-1}(a_1^k) - 2 + {|A|}/{2} )
\cdots( {|A|}/{2} )
}
\]
is the $k$-order KT-probability of $X_1^n$. (For $k=0$, $N_{n-1}(a_1^k)=n$.)
\end{definition}
%
%
\begin{remark}
The $k$-order KT-probability of the sample is equal to a mixture of the
probabilities of the sample with respect to all $k$-order Markov chains
with uniform initial distribution, where the mixture distribution over
the transition probability matrices $ \{ Q(a_{k+1}|a_1^k),
a_{k+1}\in A, a_1^k\in A^k \}$ is independent for the rows
$Q( \cdot|a_1^k)$, $a_1^k\in A^k$, and has
Dirichlet $(
\frac12,\ldots,\frac12 )$ distribution in the rows. Hence,
the KT
Markov order estimator can be regarded as a Bayes (maximum a
posteriori) estimator.
\end{remark}
%
%
\begin{remark}
The $k$-order NML and KT coding distributions are nearly optimal among
the $k$-order Markov chains, in the sense that the code lengths $\lceil
-\log P_{\mathrm{NML}, k} (X_1^n) \rceil$ and $\lceil-\log
P_{\mathrm{KT}, k} (X_1^n) \rceil$ minimize the
worst case maximum and average, respectively, redundancy for this class
(up to an additive constant in the latter case).
\end{remark}

\section{Divergence of Markov order estimators}\label{secmain}

The BIC Markov order estimator is strongly consistent~\cite{CsSh}, that
is, if the process is a Markov chain of order $k$, then $\hat
{k}_{\mathrm{BIC}} (X_1^n) = k$ eventually almost surely as $n\to
+\infty
$. ``Eventually almost surely'' means that with probability $1$, there
exists a threshold $n_0$ (depending on the infinite realization
$X_{1}^{\infty}$) such that the claim holds for all $n \ge n_0$.
Increasing the penalty term, up to $cn$, where $c>0$ is a sufficiently
small constant, does not affect the strong consistency. It is not known
whether or not the strong consistency holds for smaller penalty terms
but it is known that if the candidate orders are upper bounded by
$c\log n$, where $c>0$ is a sufficiently small constant, that is, the
estimator minimizes the PML over the orders $0\le k \le c\log n$ only,
then $\operatorname{pen}(n) = C \log\log n$ still provides the strong
consistency,
where $C>0$ is a sufficiently large constant~\cite{Ramon}.

The NML and KT Markov order estimators fail to be strongly consistent
because for i.i.d. processes with uniform distribution, they converge
to infinity at a rate $\mathcal{O}(\log n)$~\cite{CsSh}. However, if
the candidate orders are upper bounded by $\RMo(\log n)$, the strong
consistency holds true~\cite{Cs}.

If the process is of infinite memory, the BIC and KT Markov order
estimators diverge to infinity~\cite{CsT}. In this section, results on
the divergence rate of the PML, NML and KT Markov order estimators are
presented. Bounds on the probability that the estimators are greater
and less, respectively, than some order are obtained, with explicit
constants. The first implies that under mild conditions, the estimators
do not exceed the $\mathcal{O}(\log n)$ rate eventually almost surely
as $n\to+\infty$. The second bound implies that the rate $\mathcal
{O}(\log n)$ is attained eventually almost surely as $n\to+\infty$ for
the processes whose continuity rates decay in some exponential range.

At the end of the section, the notion of consistent Markov order
estimation is generalized for infinite memory processes. If the
continuity rates decay faster than exponential, the PML Markov order
estimator is shown to be consistent with the oracle-type order estimate.

The proofs use bounds on the simultaneous convergence of empirical
entropies of orders in an increasing set. These bounds are obtained for
finite sample sizes $n$ with explicit constants under mild conditions
so they are of independent interest and are also presented here.
%
%
\begin{theorem}\label{thentSMP}
For any weakly non-null and $\alpha$-summable stationary ergodic
process, for any $0<\varepsilon<1/2$
\[
\Prob\biggl( \max_{1\le k\le({\varepsilon\log n})/({4\log|A|})} \bigl
\llvert\hat{H}_{k}
\bigl(X_1^n\bigr) - H_k \bigr\rrvert>
\frac{1}{n^{1/2-\varepsilon}} \biggr) \le\exp\biggl( -\frac{c_1\varepsilon^3}{\log
n} n^{\varepsilon/2}
\biggr)
\]
and
\[
\Prob\biggl( \max_{0\le k\le({\varepsilon\log n})/({4\log|A|})} \bigl
\llvert\hat{h}_{k}
\bigl(X_1^n\bigr) - h_k \bigr\rrvert>
\frac{1}{n^{
1/2-\varepsilon}} \biggr) \le\exp\biggl( -\frac{c_2\varepsilon^3}{\log
n} n^{\varepsilon/2}
\biggr),
\]
where $c_1,c_2>0$ are constants depending only on the distribution of
the process.
\end{theorem}
\begin{pf}
The proof including the explicit expression of the constants is in
Section~\ref{secent}.
\end{pf}
%
%
\begin{remark}
The convergence of $\hat{H}_{k_n}(X_1^n)$ and $\hat{h}_{k_n}(X_1^n)$,
$k_n\to\infty$, to the
entropy rate $\bar{H}$ of the process could be investigated using
Theorem~\ref{thentSMP}. However, good estimates of the entropy rate
are known from the theory of universal codes. In particular, mixtures
of the KT distributions over all possible orders provide universal
codes in the class of all stationary ergodic processes \cite
{Ryabko,Ryabko2,Ryabko3}, therefore the corresponding code length is a
suitable estimate of the entropy rate.
\end{remark}

An application of the Borel--Cantelli lemma in Theorem~\ref{thentSMP}
yields the following asymptotic result.
%
%
\begin{corollary}
For any weakly non-null and $\alpha$-summable stationary ergodic
process, for any $0<\varepsilon<1/2$
\[
\bigl\llvert\hat{H}_{k}\bigl(X_1^n\bigr)
- H_k \bigr\rrvert\le\frac{1}{n^{
1/2-\varepsilon}} \quad\mbox{and}\quad\bigl
\llvert\hat{h}_{k}\bigl(X_1^n\bigr) -
h_k \bigr\rrvert\le\frac{1}{n^{
1/2-\varepsilon}}
\]
simultaneously for all $k\le\frac{\varepsilon\log n}{4\log|A|}$,
eventually almost surely as $n\to+\infty$.
\end{corollary}
%
%
\begin{remark}
By~\cite{GGG}, under much stronger conditions on the process, the
convergence rate of $\hat{H}_{k}(X_1^n)$ and $\hat{h}_{k}(X_1^n)$ to
$\bar{H}$ is $n^{-1/2}$
for some fixed $k=\mathcal{O}(\log n)$. Hence, the rate in Theorem
\ref
{thentSMP} cannot be improved significantly.
\end{remark}

The first divergence result of the paper is the following.
%
%
\begin{theorem}\label{thpmllarge}
For any weakly non-null and $\alpha$-summable stationary ergodic
process there exist $\lambda_1, \lambda_2>0$ depending only on the
distribution of the process, such that for the Markov order estimator
$\hat{k}_{\mathrm{IC}} (X_1^n)$
\[
\Prob\bigl( \hat{k}_{\mathrm{IC}} \bigl(X_1^n\bigr)
> k_n \bigr) \le2^{\lambda_1 + 2\log n - \lambda_2 k_n}
\]
for any sequence $k_n$, $n\in\mathbb{N}$, where IC is either the PML
with arbitrary $\operatorname{pen}(n)$ or the NML or the KT criterion.
\end{theorem}
\begin{pf}
The proof including the explicit expression of the constants is in
Section~\ref{secproof}.
\end{pf}

An application of the Borel--Cantelli lemma in Theorem \ref
{thpmllarge} yields the following asymptotic result.
%
%
\begin{corollary}\label{copmllarge}
For any weakly non-null and $\alpha$-summable stationary ergodic
process there exists a constant $C>0$ such that for the Markov order
estimator $\hat{k}_{\mathrm{IC}} (X_1^n)$
\[
\hat{k}_{\mathrm{IC}} \bigl(X_1^n\bigr) \le C \log n
\]
eventually almost surely as $n\to+\infty$, where IC is either the PML
with arbitrary $\operatorname{pen}(n)$ or the NML or the KT criterion.
\end{corollary}

The second divergence result is the following.
%
%
\begin{theorem}\label{thpmllogSMP}
For any weakly non-null stationary ergodic process with continuity
rates $\bar{\gamma}(k) \le\delta_1 2^{-\zeta_1 k}$ and $\underbar
{\gamma}(k) \ge\delta_2 2^{-\zeta_2 k}$ for some $\zeta_1, \zeta_2,
\delta_1, \delta_2 >0$ ($\zeta_2\ge\zeta_1$), if
\[
\frac{6\log|A|}{\zeta_1} \le\varepsilon< \frac12,
\]
the Markov order estimator $\hat{k}_{\mathrm{IC}} (X_1^n)$ satisfies that
\[
\Prob\biggl( \hat{k}_{\mathrm{IC}} \bigl(X_1^n
\bigr) \le\frac{1}{2\zeta_2} \biggl( \frac12 -\varepsilon\biggr) \log
n -
c_3 \biggr) \le\exp\biggl( -\frac{c_2\varepsilon^3}{\log n}
n^{\varepsilon/2}
\biggr),
\]
if $n\ge n_0$, where IC is either the PML with $\operatorname
{pen}(n)\le\mathcal
{O}(\sqrt{n})$ or the NML or the KT criterion,
and $c_2, n_0>0$, $c_3\in\mathbb{R}$ are constants depending only on
the distribution of the process and $\operatorname{pen}(n)$.
\end{theorem}
\begin{pf}
The proof including the explicit expression of the constants is in
Section~\ref{secproof}.
\end{pf}

An application of the Borel--Cantelli lemma in Theorem \ref
{thpmllogSMP} yields the following asymptotic result.
%
%
\begin{corollary}\label{corpmllog}
For any weakly non-null stationary ergodic process with continuity
rates $\bar{\gamma}(k) \le\delta_1 2^{-\zeta_1 k}$ and $\underbar
{\gamma}(k) \ge\delta_2 2^{-\zeta_2 k}$ for some $\zeta_1, \zeta_2,
\delta_1, \delta_2 >0$ with $\zeta_2\ge\zeta_1>12\log|A|$,
the Markov order estimator $\hat{k}_{\mathrm{IC}} (X_1^n)$ satisfies that
\[
\hat{k}_{\mathrm{IC}} \bigl(X_1^n\bigr) \ge
C' \log n
\]
eventually almost surely as $n\to+\infty$, where IC is either the PML
with $\operatorname{pen}(n)\le\mathcal{O}(\sqrt{n})$ or the NML or
the KT criterion,
and $C'>0$ is a constant depending only on the distribution of the process.
\end{corollary}

The section concludes with the consistency result.
%
%
\begin{definition}\label{defoPML}
For a candidate order $0\le k<n$ the oracle PML criterion is
\[
\mathrm{PML}_{o,n} ( k ) = (n-k) h_k + \bigl(|A|-1\bigr)
|A|^k \operatorname{pen}(n),
\]
and the oracle PML Markov order estimator is
\[
k_{\mathrm{PML},n} = \arg\min_{0\le k < n} \mathrm{PML}_{o,n} ( k ).
\]
\end{definition}
%
%
\begin{remark}
For Markov chains of order $k$, $k_{\mathrm{PML},n} = k$ if $n$ is
sufficiently large, with any $\operatorname{pen}(n) = \RMo(n)$.
\end{remark}
%
%
\begin{theorem}\label{thoracle}
For any weakly non-null stationary ergodic process with
\[
\frac{\log\bar{\gamma}(k)}{k} \to-\infty,\qquad k\to\infty,
\]
the PML Markov order estimator $\hat{k}_{\mathrm{PML}} (X_1^n)$ with
$\operatorname{pen}(n)=n^{\kappa}$, $\frac12<\kappa<1$, is
consistent in the sense that
\[
\frac{ \hat{k}_{\mathrm{PML}} (X_1^n) }{ k_{\mathrm{PML},n} } \to1
\]
in probability as $n\to+\infty$.
\end{theorem}
\begin{pf}
The proof is in Section~\ref{secproof}.
\end{pf}

\section{Statistical estimation of processes}\label{secappl}

In the results of this section, the divergence rate of Markov order
estimators will play a central role. The problem of statistical
estimation of stationary ergodic processes by finite memory processes
is considered, and the following distance is used. The per-letter
Hamming distance between two strings $x_1^n$ and $y_1^n$ is
\[
d_n \bigl(x_1^n, y_1^n
\bigr) = \frac1n \sum_{i=1}^n
\mathbb{I}(x_i\neq y_i) \qquad\mbox{where } \mathbb{I}(a\neq
b) = \cases{ %
1, &\quad if $a\neq b$,
\cr
0, &\quad if $a=b$,}
\]
and the \textit{$\bar{d}$-distance} between two random sequences $X_1^n$
and $Y_1^n$ is defined by
\[
\bar{d} \bigl(X_1^n, Y_1^n
\bigr) = \min_{\mathbb{P}} \mathbb{E}_{\mathbb
{P}} d_n \bigl(
\tilde{X}_1^n, \tilde{Y}_1^n
\bigr),
\]
where the minimum is taken over all the joint distributions $\mathbb
{P}$ of $\tilde{X}_1^n$ and $\tilde{Y}_1^n$ whose marginals are equal
to the distributions of $X_1^n$ and $Y_1^n$.

The process $X$ is estimated by a Markov chain of order $k=k_n$ from
the sample in the following way.
%
%
\begin{definition}
The \textup{empirical $k$-order Markov estimator} of a process $X$ based
on the sample $X_1^n$ is the stationary Markov chain, denoted by $\hat
{X}[k]$, of order $k$ with transition probability matrix $ \{
\hat
{P} (a_{k+1}| a_1^k), a_{k+1}\in A, a_1^k\in A^k \}$. If the
initial distribution of a stationary Markov chain with these transition
probabilities is not unique, then any of these initial distributions
can be taken.
\end{definition}

In the previous section, weakly non-nullness is assumed for the
process. In this section the process $X$ is assumed to be \textit
{non-null}, that is,
\[
\pinf= \min_{a\in A} \inf_{x_{-\infty}^{-1} \in A^{\infty}} P\bigl(a |
x_{-\infty}^{-1}
\bigr) >0.
\]
%
%
\begin{remark}\label{rempmllarge}
For any non-null stationary ergodic process, $P(a_1^k)\le(1-\pinf)^k$
for any $a_1^k\in A^k$. Hence, Theorem~\ref{thpmllarge} holds with
$\lambda_1=0$ and $\lambda_2=|{\log}(1-\pinf)|$, see the proof of the theorem.
\end{remark}

The assumption of non-nullness allows us to use the following quantity
instead of $\underbar{\gamma}(k)$.
The \textit{restricted continuity rate} of the process $X$ is
\[
\bar{\gamma}(k|m) = \max_{x_{-m}^{-1} \in A^m} \sum_{a\in A}
\bigl\llvert P\bigl(a | x_{-k}^{-1} \bigr) - P\bigl(a |
x_{-m}^{-1} \bigr) \bigr\rrvert,\qquad k<m.
\]
Similarly to Remark~\ref{remgammaeq}, note that the above definition
is equivalent to
\[
\bar{\gamma}(k|m) = \max_{ k<i\le m } \max_{x_{-i}^{-1} \in A^i} \sum
_{a\in A} \bigl\llvert P\bigl(a | x_{-k}^{-1}
\bigr) - P\bigl(a | x_{-i}^{-1} \bigr) \bigr\rrvert.
\]
Hence, $\lim_{m\to+\infty} \bar{\gamma}(k|m) = \bar{\gamma}(k)$
for any
fixed $k$. We say that the process $X$ has \textit{uniformly convergent
restricted continuity rate} with parameters $\theta_1$, $\theta_2$,
$k_{\theta}$ if
\[
\bar{\gamma}(k)^{\theta_1} \le\bar{\gamma}\bigl( k | \lceil
\theta_2 k \rceil\bigr) \qquad\mbox{if } k\ge k_{\theta}
\mbox{, for some } \theta_1\ge1,\theta_2>1. 
\]

The order $k$ of the empirical Markov estimator $\hat{X}[k]$ is
estimated from the sample, using the PML criterion. The estimated order
needs to be bounded to guarantee an accurate assessment of the memory
decay of the process.
%
%
\begin{definition}\label{defICr}
For an information criterion IC, the Markov order estimator bounded by
$r_n<n$, $r_n\in\mathbb{N}$, is
\[
\hat{k}_{\mathrm{IC}} \bigl(X_1^n | r_n
\bigr) = \arg\min_{0\le k \le r_n} \mathrm{IC}_{X_1^n} (k).
\]
\end{definition}

The optimal order can be smaller than the upper bound if the memory
decay of the process is sufficiently fast. Define
\[
K_n \bigl( r_n, \bar{\gamma}, f(n) \bigr) = \min\bigl\{
\lfloor r_n \rfloor, k\ge0\dvtx\bar{\gamma}(k) < f(n) \bigr\} ,
\]
where $f(n)\searrow0$ and $r_n\nearrow\infty$. Since $\bar{\gamma
}$ is
a decreasing function, $K_n$ increases in $n$ but does not exceed
$r_n$. It is less than $r_n$ if $\bar{\gamma}$ vanishes sufficiently
fast, and then the faster $\bar{\gamma}$ vanishes, the slower $K_n$ increases.

The process estimation result of the paper is the following.
%
%
\begin{theorem}\label{thqminSMP}
For any non-null stationary ergodic process with summable continuity
rate and uniformly convergent restricted continuity rate with
parameters $\theta_1$, $\theta_2$, $k_{\theta}$, and for any $\mu_n>0$,
the empirical Markov estimator of the process with the order estimated
by the bounded PML Markov order estimator $\hat{k}_n = \hat
{k}_{\mathrm
{PML}} (X_1^n | \eta\log n)$, $\eta>0$, with $\frac12 \log n \le
\operatorname{pen}(n)\le\mathcal{O}(\sqrt{n})$ satisfies
\begin{eqnarray*}
&&\Prob\biggl( \bar{d} \bigl( X_1^n, \hat{X} [
\hat{k}_n ]_1^n \bigr) > \frac{\beta_2}{\pinf^{2}}
\max\biggl\{ \bar{\gamma} \biggl( \biggl\lfloor\frac{\eta}{\theta_2}
\log n \biggr
\rfloor\biggr), n^{-( 1-4\eta\log({|A|^4}/{\pinf})
)/({4\theta_1})
} \biggr\} + \frac{1}{n^{1/2-\mu_n}} \biggr)
\\
&&\quad\le\exp\bigl( -c_4 4^{ \mu_n\log n - |\log\pinf| (
K_n ( \eta\log n, \bar{\gamma}, {c}\operatorname
{pen}(n)/{n} ) +
{\log\log n}/{\log|A|} ) } \bigr)
\\
&&\qquad{} + \exp\biggl( -\frac{c_5 \eta^3}{\log n} n^{\eta2\log|A|}
\biggr) +
2^{-s_n \operatorname{pen}(n)},
\end{eqnarray*}
if $n\ge n_0$, where $c>0$ is an arbitrary constant, $s_n\to\infty$ and
$\beta_2, c_4, c_5, n_0>0$ are constants depending only on the
distribution of the process.
\end{theorem}
\begin{pf}
The proof including the explicit expression of the constants is in
Section~\ref{secproofappl}.
\end{pf}
%
%
\begin{remark}
If the process $X$ is a Markov chain of order $k$, then the restricted
continuity rate is uniformly convergent with parameters $\theta_1=1$,
$\theta_2>1$ arbitrary (arbitrarily close to $1$), $k_{\theta}=k+1$,
and if $n$ is sufficiently large, $K_n=k$ and
\[
\max\biggl\{ \bar{\gamma} \biggl( \biggl\lfloor\frac{\eta
}{\theta_2} \log n \biggr
\rfloor\biggr), n^{- (
1-4\eta
\log({|A|^4}/{\pinf}) )/({4\theta_1}) } \biggr\} = n^{- ( 1-4\eta
\log({|A|^4}/{\pinf})
)/({4\theta_1}) } .
\]
\end{remark}

An application of the Borel--Cantelli lemma in Theorem~\ref{thqminSMP}
yields the following asymptotic result.
%
%
\begin{corollary}\label{cormain}
For any non-null stationary ergodic process with summable continuity
rate and uniformly convergent restricted continuity rate with
parameters $\theta_1$, $\theta_2$, $k_{\theta}$,
the empirical Markov estimator of the process with the order estimated
by the bounded PML Markov order estimator $\hat{k}_n = \hat
{k}_{\mathrm
{PML}} (X_1^n | r_n)$ with $\frac12 \log n \le\operatorname
{pen}(n)\le\mathcal
{O}(\sqrt{n})$ and
\[
\frac{5\log\log n}{2\log|A|} \le r_n \le\RMo(\log n)
\]
satisfies
\begin{eqnarray*}
\bar{d} \bigl( X_1^n, \hat{X} [ \hat{k}_n
]_1^n \bigr) &\le&\frac{\beta_2}{\pinf^{2}} \max\biggl\{ \bar{
\gamma} \biggl( \biggl\lfloor\frac{r_n}{\theta_2} \biggr\rfloor\biggr),
n^{-
{1}/({4\theta_1}) } \biggr\}\\
&&{} + \frac{ (\log n)^{c_6} }{ \sqrt{n} } 2^{
|\log\pinf| K_n ( r_n,
\bar{\gamma}, {c}\operatorname{pen}(n)/{n} ) }
\end{eqnarray*}
eventually almost surely as $n\to+\infty$,
where $c>0$ is an arbitrary constant, and $\beta_2, c_6>0$ are
constants depending only on the distribution of the process.
\end{corollary}
%
%
\begin{remark}
If the memory decay of the process is slow, the first term in the bound
in Corollary~\ref{cormain}, the bias, is essentially $\bar{\gamma} (
\lfloor r_n/\theta_2 \rfloor)$, and the second term, the
variance, is maximal. If the memory decay is sufficiently fast, then
the rate of the estimated order $\hat{k}_n$ and the rate of $K_n$ are
smaller, therefore the variance term is smaller, while the bias term is
smaller as well. The result, however, shows the optimality of the PML
Markov order estimator in the sense that it selects an order which is
small enough to allow the variance to decrease but large enough to keep
the bias below a polynomial threshold.
\end{remark}

\section{Empirical entropies}\label{secent}

In this section, we consider the problem of simultaneous convergence of
empirical entropies of orders in an increasing set, and prove the
following theorem that formulates Theorem~\ref{thentSMP} with explicit
constants.
%
%
\begin{theorem}\label{thent}
For any weakly non-null and $\alpha$-summable stationary ergodic
process, for any $0<\varepsilon<1/2$
\begin{eqnarray*}
&&
\Prob\biggl( \max_{1\le k\le({\varepsilon\log n})/({4\log|A|})} \bigl
\llvert\hat{H}_{k}
\bigl(X_1^n\bigr) - H_k \bigr\rrvert>
\frac{1}{n^{
1/2-\varepsilon}} \biggr) \\
&&\quad\le6 \RMe ^{1/\RMe } \exp\biggl( -\frac{7\alpha
_0\varepsilon^3}{32\RMe (\alpha
+\alpha_0)}
\frac{ n^{\varepsilon/2} }{ \log n } + \frac
{\varepsilon
}{4} \log n \biggr)
\end{eqnarray*}
and
\begin{eqnarray*}
&&
\Prob\biggl( \max_{0\le k\le({\varepsilon\log n})/({4\log|A|})} \bigl
\llvert\hat{h}_{k}
\bigl(X_1^n\bigr) - h_k \bigr\rrvert>
\frac{1}{n^{
1/2-\varepsilon}} \biggr) \\
&&\quad\le12 \RMe ^{1/\RMe } \exp\biggl( -\frac{7\alpha
_0\varepsilon^3}{256\RMe (\alpha+\alpha_0)}
\frac{ n^{\varepsilon/2} }{ \log n } + \frac{\varepsilon}{4} \log n
\biggr) .
\end{eqnarray*}
\end{theorem}

First, we show the following bounds.
%
%
\begin{proposition}\label{prent}
For any weakly non-null and $\alpha$-summable stationary ergodic
process, for any $1\le m\le n$ and $u,\nu>0$,
\begin{eqnarray*}
&& \Prob\Bigl( \max_{1\le k\le m} \bigl\llvert\hat{H}_{k}
\bigl(X_1^n\bigr) - H_k \bigr\rrvert> u
\Bigr)
\\
&&\quad\le6 \RMe ^{1/\RMe } |A|^m \exp\biggl( \frac{\alpha_0}{8\RMe (\alpha
+\alpha_0)}
\frac{-(n-m+1)u^{2(1+\nu)}}{ m |A|^{2m} }
\\
&&\qquad\hspace*{60.5pt}{} \times\min{}^2 \biggl\{ \biggl( \frac{\RMe }{2(1+\nu^{-1})}
\biggr)^{1+\nu}, \frac{u^{-\nu}\log \RMe }{2m\log|A|}, \frac{u^{-\nu}}{\RMe }
\biggr\}\biggr)
\end{eqnarray*}
and
\begin{eqnarray*}
&& \Prob\Bigl( \max_{0\le k\le m-1} \bigl\llvert\hat{h}_{k}
\bigl(X_1^n\bigr) - h_k \bigr\rrvert> u
\Bigr)
\\
&&\quad\le12 \RMe ^{1/\RMe } |A|^m \exp\biggl( \frac{\alpha_0}{8\RMe (\alpha
+\alpha_0)}
\frac{-(n-m+1) (u/2)^{2(1+\nu)}}{ m |A|^{2m} }
\\
&&\qquad\hspace*{64.5pt}{} \times \min{}^2 \biggl\{ \biggl( \frac{\RMe }{2(1+\nu^{-1})}
\biggr)^{1+\nu}, \frac{(u/2)^{-\nu}\log \RMe }{2m\log|A|}, \frac{(u/2)^{-\nu
}}{\RMe } \biggr\} \biggr).
\end{eqnarray*}
\end{proposition}
\begin{pf}
Fix $1\le k\le m$. Applying Lemma~\ref{lemsch} in the \hyperref
[app]{Appendix} to the
distributions $P_k = \{ P (a_1^k), a_1^k\in A^k \}$ and
$\hat{P}_k = \{ \hat{P} (a_1^k), a_1^k\in A^k \}$,
%
%
\begin{equation}
\label{eqsch} \bigl\llvert\hat{H}_{k}\bigl(X_1^n
\bigr) - H_k \bigr\rrvert\le\frac{1}{\log \RMe } \bigl[ k\log|A| - \log
d_{\TV} ( \hat{P}_k, P_k ) \bigr]
d_{\TV} ( \hat{P}_k, P_k ),
\end{equation}
if $d_{\TV} ( \hat{P}_k, P_k ) \le1/\RMe $. For any $\nu>0$, the right
of (\ref{eqsch}) can be written as
%
%
\begin{eqnarray}
\label{eqsch2} && \frac{k \log|A|}{\log \RMe } d_{\TV} ( \hat{P}_k,
P_k ) \nonumber\\
&&\qquad{}+ \frac{1+\nu}{\nu\log \RMe } d_{\TV}^{{1}/({1+\nu})} (
\hat{P}_k , P_k ) \bigl[ - d_{\TV}^{{\nu}/({1+\nu})}
( \hat{P}_k, P_k ) \log d_{\TV}^{{\nu}/({1+\nu})}
( \hat{P}_k, P_k ) \bigr]
\\
&&\quad\le\frac{k \log|A|}{\log \RMe } d_{\TV} ( \hat{P}_k,
P_k ) + \frac{1}{\RMe } \frac{1+\nu}{\nu} d_{\TV}^{{1}/({1+\nu})}
( \hat{P}_k, P_k ),\nonumber
\end{eqnarray}
where we used the bound $- x\log x \le \RMe ^{-1} \log \RMe $, $x\ge0$.

By~\cite{GL}, for any string $a_1^k\in A^k$ and $t>0$,
%
%
\begin{equation}
\label{eqGL} \Prob\bigl( \bigl\llvert N_n\bigl(a_1^k
\bigr) - (n-k+1) P\bigl(a_1^k\bigr) \bigr\rrvert> t
\bigr) \le \RMe ^{1/\RMe } \exp\biggl( \frac{-c_{\alpha} t^2}{ k(n-k+1) }
\biggr),
\end{equation}
where
\[
c_{\alpha}= \frac{\alpha_0}{8\RMe (\alpha+\alpha_0)}
\]
is positive for any weakly non-null and $\alpha$-summable stationary
ergodic process.
(\ref{eqGL}) implies that
%
%
\begin{eqnarray}
\label{eqGL2} \Prob\bigl( d_{\TV} ( \hat{P}_k,
P_k ) > t \bigr) &\le&\Prob\biggl( \max_{a_1^k\in A^k} \bigl\llvert
\hat{P} \bigl(a_1^k\bigr) - P\bigl(a_1^k
\bigr) \bigr\rrvert> \frac{t}{|A|^k} \biggr)
\nonumber\\[-8pt]\\[-8pt]
&\le& \RMe ^{1/\RMe } |A|^k \exp\biggl( \frac{-c_{\alpha}(n-k+1)
t^2}{ k |A|^{2k} }
\biggr).\nonumber
\end{eqnarray}
Applying (\ref{eqGL2}) to (\ref{eqsch2}),
\begin{eqnarray*}
&&\Prob\bigl( \bigl\llvert\hat{H}_{k}
\bigl(X_1^n\bigr) - H_k \bigr\rrvert> u
\bigr)
\\
&&\quad\le\Prob\biggl( \frac{k \log|A|}{\log \RMe } d_{\TV} (
\hat{P}_k , P_k ) + \frac{1}{\RMe } \frac{1+\nu}{\nu}
d_{\TV}^{{1}/({1+\nu})} ( \hat{P}_k, P_k ) > u
\biggr)
\\
&&\qquad{} + \Prob\bigl( d_{\TV} ( \hat{P}_k,
P_k ) > 1/\RMe  \bigr)
\\
&&\quad\le\Prob\biggl( d_{\TV} ( \hat{P}_k,
P_k ) > \frac{u\log \RMe }{2k\log|A|} \biggr) + \Prob\biggl(
d_{\TV}^{{1}/({1+\nu})} ( \hat{P}_k, P_k ) >
\frac
{\nu \RMe  u}{2(1+\nu)} \biggr)
\\
&&\qquad{} + \Prob\bigl( d_{\TV} ( \hat{P}_k,
P_k ) > 1/\RMe  \bigr)
\\
&&\quad\le3 \RMe ^{1/\RMe } |A|^k \exp\biggl( \frac{-c_{\alpha
}(n-k+1)u^{2(1+\nu
)}}{ k |A|^{2k} }\\
&&\qquad\hspace*{59.5pt}{}\times
\min{}^2 \biggl\{ \biggl( \frac{\RMe }{2(1+\nu^{-1})} \biggr)^{1+\nu} ,
\frac{u^{-\nu}\log \RMe }{2k\log|A|}, \frac{u^{-\nu}}{\RMe } \biggr\} \biggr).
\end{eqnarray*}
This completes the proof of the first claimed bound as
\begin{eqnarray*}
&&\Prob\Bigl( \max_{1\le k\le m} \bigl\llvert\hat{H}_{k}
\bigl(X_1^n\bigr) - H_k \bigr\rrvert> u
\Bigr)
\\
&&\quad\le\sum_{1\le k\le m} \Prob\bigl( \bigl\llvert
\hat{H}_{k}\bigl(X_1^n\bigr) -
H_k \bigr\rrvert> u \bigr)
\\
&&\quad\le3 \RMe ^{1/\RMe } \biggl( \sum_{1\le k\le m}|A|^k
\biggr)
\\
&&\qquad\hspace*{0pt}{}\times
\exp\biggl( \frac{-c_{\alpha}(n-m+1)u^{2(1+\nu)}}{ m
|A|^{2m} } \min{}^2 \biggl\{
\biggl( \frac{\RMe }{2(1+\nu^{-1})} \biggr)^{1+\nu}, \frac{u^{-\nu}\log
\RMe }{2m\log|A|},
\frac{u^{-\nu}}{\RMe } \biggr\} \biggr).
\end{eqnarray*}
The second claimed bound follows using $\hat{h}_{0}(X_1^n) - h_0 =
\hat{H}_{1}(X_1^n) - H_1$ and
\[
\bigl\llvert\hat{h}_{k}\bigl(X_1^n\bigr)
- h_k \bigr\rrvert\le\bigl\llvert\hat{H}_{k+1}
\bigl(X_1^n\bigr) - H_{k+1} \bigr\rrvert+
\bigl\llvert\hat{H}_{k}\bigl(X_1^n\bigr) -
H_k \bigr\rrvert,\qquad k\ge1,
\]
as
\begin{eqnarray*}
&& \Prob\Bigl( \max_{0\le k\le m-1} \bigl\llvert\hat{h}_{k}
\bigl(X_1^n\bigr) - h_k \bigr\rrvert> u
\Bigr)
\\
&&\quad\le\Prob\biggl( \max_{1\le k\le m} \bigl\llvert\hat{H}_{k}
\bigl(X_1^n\bigr) - H_k \bigr\rrvert>
\frac{u}{2} \biggr) + \Prob\biggl( \max_{1\le k\le m-1} \bigl\llvert
\hat{H}_{k}\bigl(X_1^n\bigr) -
H_k \bigr\rrvert> \frac{u}{2} \biggr)
\\
&&\quad\le2\Prob\biggl( \max_{1\le k\le m} \bigl\llvert\hat{H}_{k}
\bigl(X_1^n\bigr) - H_k \bigr\rrvert>
\frac{u}{2} \biggr) .
\end{eqnarray*}
\upqed\end{pf}

Now, the theorem follows from the proposition with special settings.
\begin{pf*}{Proof of Theorem~\ref{thent}}
We use Proposition~\ref{prent} setting $u = n^{-1/2+\varepsilon}$,
$\nu= \varepsilon$, and $m = \lfloor(\varepsilon\log n)/(4\log|A|)
\rfloor$. Then, in the exponent of the first inequality of the proposition,
\begin{eqnarray*}
\frac{ u^{2(1+\nu)} }{ |A|^{2m} } &>& n^{\varepsilon/2 -1
+2\varepsilon
^2},
\\
\frac{n-m+1}{m} &>& n \frac{7}{\log n},
\\
\min\biggl\{ \biggl( \frac{\RMe }{2(1+\nu^{-1})} \biggr)^{1+\nu},
\frac
{u^{-\nu}\log \RMe }{2m\log|A|}, \frac{u^{-\nu}}{\RMe } \biggr\} &>& n^{
-\varepsilon^2 } \biggl(
\frac{2\varepsilon}{3} \biggr)^{3/2} > n^{ -\varepsilon^2 } \frac
{\varepsilon^{3/2}}{2},
\end{eqnarray*}
where we used that $0<\varepsilon<1/2$. This gives the lower bound
\[
-\frac{7\alpha_0\varepsilon^3}{32\RMe (\alpha+\alpha_0)} \frac{
n^{\varepsilon/2} }{ \log n }
\]
on the exponent and completes the proof of the first claimed bound. The
second claimed bound follows similarly from the second inequality of
the proposition with the same settings.
\end{pf*}

\section{Divergence bounds proofs}\label{secproof}

In this section, we consider the divergence of the PML, NML and KT
Markov order estimators and prove Theorems~\ref{thpmllarge}, \ref
{thpmllogSMP} and~\ref{thoracle}.
\begin{pf*}{Proof of Theorem~\ref{thpmllarge}}
By~\cite{GL}, any weakly non-null and $\alpha$-summable process is
$\phi
$-mixing with a coefficient related to $\alpha_0>0$ and $\alpha
<+\infty
$. Namely, there exists a sequence $\rho_i$, $i\in\mathbb{N}$, satisfying
\[
\sum_{i=0}^{\infty} \rho_i \le1+
\frac{2\alpha}{\alpha_0} ,
\]
such that for each $k$, $m$, $l$ and each $a_1^k\in A^k$, $b_1^m\in
A^m$, with $P(b_1^m)>0$,
\[
\bigl\llvert\Prob\bigl( X_{m+l+1}^{m+l+k} = a_1^k
\mid X_1^m = b_1^m \bigr) - P
\bigl(a_1^k\bigr) \bigr\rrvert\le\sum
_{i=l}^{l+k-1} \rho_i .
\]
This implies that for any $d\ge1$
\begin{eqnarray*}
&& \Prob\bigl( X_{m+l+1}^{m+l+k} = a_1^k
\mid X_1^m = b_1^m \bigr)
\\
&&\quad\le\Prob\bigl( X_{m+l+id} = a_{id}, 1\le i\le\lfloor
k/d \rfloor\mid X_1^m = b_1^m
\bigr)
\\
&&\quad= \prod_{i=1}^{ \lfloor k/d \rfloor} \Prob\bigl(
X_{m+l+id} = a_{id} \mid X_{m+l+jd} = a_{jd},
1\le j<i, X_1^m = b_1^m
\bigr)
\\
&&\quad\le\prod_{i=1}^{ \lfloor k/d \rfloor} \bigl(
P(a_{id}) + \rho_{d-1} \bigr)
\\
&&\quad\le\Bigl( \max_{a\in A} P(a) + \rho_{d-1}
\Bigr)^{ \lfloor k/d
\rfloor}.
\end{eqnarray*}
Since $\max_{a\in A} P(a)<1$ and $\rho_d\to0$, $\max_{a\in A} P(a) +
\rho_{d-1} <1$ for sufficiently large $d$. Then
\[
\max_{l,a_1^k,b_1^m} \Prob\bigl( X_{m+l+1}^{m+l+k} =
a_1^k \mid X_1^m =
b_1^m \bigr) \le2^{\lambda_1 - \lambda_2 k }
\]
holds with $\lambda_1 = -\log( \max_{a\in A} P(a) + \rho_{d-1}
) >0$ and $\lambda_2 = -\log( \max_{a\in A}
P(a) +
\rho_{d-1} )^{1/d} >0$. Thus, for any $k$,
%
%
\begin{eqnarray}
\label{eqexp} && \Prob\bigl( N_n\bigl(a_1^k
\bigr) \ge2 \mbox{ for some } a_1^k \bigr)
\nonumber
\\
&&\quad= \Prob\bigl( X_i^{i+k-1} = X_j^{j+k-1}
\mbox{ for some } 1\le i<j\le n-k+1 \bigr)
\nonumber
\\
&&\quad\le\sum_{1\le i<j\le n-k+1} \Prob\bigl(
X_i^{i+k-1} = X_j^{j+k-1} \bigr)
\\
&&\quad= \sum_{1\le i<j\le n-k+1} \mathbb{E} \bigl\{ \Prob\bigl(
X_j^{j+k-1} = X_i^{i+k-1} |
X_1^{j-1} \bigr) \bigr\}
\nonumber
\\
&&\quad\le n^2 2^{\lambda_1 - \lambda_2 k }.\nonumber
\end{eqnarray}

For any information criterion IC, we can write
%
%
\begin{eqnarray}
\label{eqsplit} && \bigl\{ \hat{k}_{\mathrm{IC}} \bigl(X_1^n
\bigr) > k_n \bigr\}
\nonumber
\\
&&\quad\subseteq\bigl\{ \mathrm{IC}_{X_1^n} ( m ) < \mathrm
{IC}_{X_1^n} ( k_n ) \mbox{ for some } m>k_n
\bigr\}
\nonumber\\[-8pt]\\[-8pt]
&&\quad\subseteq\bigl\{ \mathrm{IC}_{X_1^n} ( m ) < \mathrm
{IC}_{X_1^n} ( k_n ) \mbox{ for some } m>k_n
\bigr\} \cap\bigl\{ N_n\bigl(a_1^{k_n}\bigr)
\le1 \mbox{ for all } a_1^{k_n} \bigr\}
\nonumber
\\
&&\qquad{} \cup\bigl\{ N_n\bigl(a_1^{k_n}
\bigr) \ge2 \mbox{ for some } a_1^{k_n} \bigr\}.\nonumber
\end{eqnarray}
Here, $N_n(a_1^{k_n}) \le1$ for all $a_1^{k_n}\in A^{k_n}$ implies that
$N_n(a_1^m) \le1$ for all $a_1^m\in A^m$ for all $m\ge k_n$, which
further implies that for all $m>k_n$ (i) $\hat{h}_{m}(X_1^n)=0$ and therefore
$\mathrm{PML}_{X_1^n} ( m ) = (|A|-1) |A|^m \operatorname
{pen}(n)$ and $\mathrm{NML}_{X_1^n} ( m ) = \Sigma_{n,m}$ and
(ii) $\mathrm{KT}_{X_1^n} ( m ) = |A|^{-n}$. Then all
the three information criteria do
not depend on the sample and are non-decreasing in $m$. Hence, in (\ref
{eqsplit})
\[
\bigl\{ \mathrm{IC}_{X_1^n} ( m ) < \mathrm{IC}_{X_1^n} (
k_n ) \mbox{ for some } m>k_n \bigr\} \cap\bigl\{
N_n\bigl(a_1^{k_n}\bigr) \le1 \mbox{ for all }
a_1^{k_n} \bigr\}
\]
is an empty set. Thus, (\ref{eqsplit}) gives
\[
\Prob\bigl( \hat{k}_{\mathrm{IC}} \bigl(X_1^n\bigr)
> k_n \bigr) \le\Prob\bigl( N_n\bigl(a_1^{k_n}
\bigr) \ge2 \mbox{ for some } a_1^{k_n} \bigr)
\]
and using (\ref{eqexp}) completes the proof.
\end{pf*}

To prove Theorem~\ref{thpmllogSMP}, first we show the
following bounds.
%
%
\begin{proposition}\label{thpmlh}
For any weakly non-null and $\alpha$-summable stationary ergodic
process with $h_k-\bar{H} \le\delta2^{-\zeta k}$ for some $\delta
,\zeta>0$, if
\[
\frac{4\log|A|}{\zeta} \le\varepsilon<\frac12 ,
\]

\begin{longlist}[(ii)]
\item[(i)]
the PML Markov order estimator $\hat{k}_{\mathrm{PML}} (X_1^n)$
satisfies that
\[
\Prob\bigl( \hat{k}_{\mathrm{PML}} \bigl(X_1^n\bigr)
< k_n \bigr) \le12 \RMe ^{1/\RMe } \exp\biggl( -\frac{7\alpha_0\varepsilon
^3}{256\RMe (\alpha+\alpha_0)}
\frac{ n^{\varepsilon/2} }{ \log n } + \frac{\varepsilon}{4} \log n
\biggr),
\]
if $n\ge(\delta2^{\zeta})^2$, where
\[
k_n = \min\biggl\{ k\ge0\dvtx h_k - \bar{H} <
\frac{ 4\max(\sqrt{n},(|A|-1) \operatorname{pen}(n)) }{
n^{1-\varepsilon} } \biggr\} ;
\]

\item[(ii)] the Markov order estimator $\hat{k}_{\mathrm{IC}} (X_1^n)$,
where IC is either NML or KT, satisfies~that
\[
\Prob\bigl( \hat{k}_{\mathrm{IC}} \bigl(X_1^n\bigr)
< k_n \bigr) \le12 \RMe ^{1/\RMe } \exp\biggl( -\frac{7\alpha_0\varepsilon
^3}{256\RMe (\alpha+\alpha_0)}
\frac{ n^{\varepsilon/2} }{ \log n } + \frac{\varepsilon}{4} \log n
\biggr),
\]
if $n\ge\max^2 \{ \sqrt{24}(\log^2 \RMe )(|A|-1)^2, 2C_{\mathrm{KT}},
\delta2^{\zeta} \}$, where
\[
k_n = \min\biggl\{ k\ge0\dvtx h_k - \bar{H} <
\frac{4}{n^{
1/2-\varepsilon}} \biggr\} .
\]
\end{longlist}
\end{proposition}
%
%
\begin{remark}
For Markov chains of order $k$, in Proposition~\ref{thpmlh} $k_n=k$
if $n$ is sufficiently large.
\end{remark}
\begin{pf*}{Proof of Proposition~\ref{thpmlh}}
Let $0<\varepsilon<1/2$ be arbitrary and
%
%
\begin{equation}
\label{eqBndef} B_n \biggl( { \frac{\varepsilon\log n}{4\log|A|} }
\biggr) = \biggl\{
\max_{0\le k\le({\varepsilon\log n})/({4\log|A|})} \bigl\llvert\hat
{h}_{k}\bigl(X_1^n
\bigr) - h_k \bigr\rrvert\le\frac{1}{n^{
1/2-\varepsilon}} \biggr\}.
\end{equation}
For any information criterion IC, we can write for any $k_n\le\frac
{\varepsilon\log n}{4\log|A|}$
%
%
\begin{eqnarray}
\label{eqsplit2} && \bigl\{ \hat{k}_{\mathrm{IC}} \bigl(X_1^n
\bigr) < k_n \bigr\}
\nonumber
\\
&&\quad\subseteq\biggl\{ \mathrm{IC}_{X_1^n} ( m ) \le\mathrm
{IC}_{X_1^n} \biggl( \biggl\lfloor{ \frac{\varepsilon\log
n}{4\log|A|} } \biggr\rfloor
\biggr) \mbox{ for some } m<k_n \biggr\}
\nonumber\\[-8pt]\\[-8pt]
&&\quad\subseteq\biggl( \biggl\{ \mathrm{IC}_{X_1^n} ( m ) \le
\mathrm{IC}_{X_1^n} \biggl( \biggl\lfloor{ \frac
{\varepsilon\log n}{4\log|A|} } \biggr
\rfloor\biggr) \mbox{ for some } m<k_n \biggr\} \cap B_n
\biggl( { \frac{\varepsilon\log n}{4\log|A|} } \biggr) \biggr)\nonumber\\
&&\qquad{} \cup
\overline{B_n
\biggl( { \frac
{\varepsilon
\log n}{4\log|A|} } \biggr)}.\nonumber
\end{eqnarray}

(i) If $\mathrm{IC} = \mathrm{PML}$, by the definition of
the PML information criterion, see Definition~\ref{defPML},
%
%
\begin{eqnarray}
\label{eqsplitpml} && \biggl\{ \mathrm{PML}_{X_1^n} ( m ) \le\mathrm
{PML}_{X_1^n} \biggl( \biggl\lfloor{ \frac{\varepsilon
\log n}{4\log|A|} } \biggr\rfloor
\biggr) \mbox{ for some } m<k_n \biggr\} \cap B_n \biggl(
{ \frac{\varepsilon\log n}{4\log|A|} } \biggr)
\nonumber
\\
&&\quad\subseteq\biggl\{ (n-m)\hat{h}_{m}\bigl(X_1^n
\bigr) - \biggl( n- \biggl\lfloor{ \frac{\varepsilon\log n}{4\log|A|}
} \biggr\rfloor\biggr) \hat
{h}_{ \lfloor({\varepsilon\log n})/({4\log|A|})
\rfloor}\bigl(X_1^n\bigr)
\nonumber
\\[-2pt]
&&\hspace*{17pt}\quad\le\bigl(|A|-1\bigr) \bigl( |A|^{ \lfloor{
({\varepsilon
\log n})/({4\log|A|}) } \rfloor} - |A|^m \bigr)
\operatorname{pen}(n) \mbox{ for some } m<k_n \biggr\} \cap
B_n \biggl( { \frac{\varepsilon\log n}{4\log|A|} } \biggr)
\nonumber
\\[-2pt]
&&\quad\subseteq\biggl\{ \hat{h}_{m}\bigl(X_1^n
\bigr) - \hat{h}_{ \lfloor
({\varepsilon\log n})/({4\log|A|}) \rfloor}\bigl(X_1^n\bigr)
\nonumber\\[-9pt]\\[-9pt]
&&\hspace*{25.5pt}\le\bigl(|A|-1\bigr) |A|^{ \lfloor{
({\varepsilon\log n})/({4\log|A|}) } \rfloor} \frac
{\operatorname{pen}(n)}{n-
\lfloor(\varepsilon\log n)/(4\log|A|) \rfloor}
\mbox{ for some } m<k_n \biggr\} \nonumber\\[-2pt]
&&\qquad{}\cap B_n \biggl( {
\frac{\varepsilon\log n}{4\log|A|} } \biggr)
\nonumber
\\[-2pt]
&&\quad\subseteq\biggl\{ h_m - h_{ \lfloor({\varepsilon\log
n})/({4\log
|A|}) \rfloor} \le
\frac{(|A|-1) |A|^{({\varepsilon\log
n})/({4\log|A|})} \operatorname{pen}(n)}{n- (\varepsilon\log n)/(4\log
|A|) } + \frac{2}{n^{1/2-\varepsilon}}\nonumber\\[-2pt]
&&\hspace*{23.5pt}\mbox{ for some } m<k_n \biggr\}.\nonumber
\end{eqnarray}
Since for any $0<\varepsilon<1/2$
%
%
\begin{equation}
\label{eqatom} \frac{ |A|^{({\varepsilon\log n})/({4\log|A|})} }{n-
(\varepsilon
\log
n)/(4\log|A|) } < \frac{1}{ n^{1-\varepsilon} },
\end{equation}
we have
%
%
\begin{equation}
\label{eqhpmlub} \frac{(|A|-1) |A|^{({\varepsilon\log n})/({4\log|A|})}
\operatorname{pen}(n)}{n-
(\varepsilon\log n)/(4\log|A|) } + \frac{2}{n^{1/2-\varepsilon}} <
\frac{ 3\max(\sqrt{n},(|A|-1) \operatorname{pen}(n)) }{
n^{1-\varepsilon} }.
\end{equation}
Now, let $\varepsilon$ and $k_n$ be as in the claim of the proposition.
Using the conditions $h_k-\bar{H} \le\delta2^{-\zeta k}$ and
$\varepsilon\ge(4\log|A|)/\zeta$,
%
%
\begin{equation}
\label{eqdhpmlub} h_{ \lfloor({\varepsilon\log n})/({4\log|A|}) \rfloor
} - \bar{H} \le\delta\exp\biggl\{ -\zeta
\biggl( \frac{\varepsilon
\log
n}{4\log|A|} -1 \biggr) \biggr\} \le\frac{1}{\sqrt{n}} \qquad
\mbox{if } n\ge\bigl(\delta2^{\zeta}\bigr)^2.
\end{equation}
Thus, if $n\ge(\delta2^{\zeta})^2$, it follows that $k_n \le\frac
{\varepsilon\log n}{4\log|A|}$, and for any $m<k_n$
%
%
\begin{eqnarray}
\label{eqhpmllb} h_m - h_{ \lfloor({\varepsilon\log n})/({4\log|A|})
\rfloor} &\ge&(h_{k_n-1} -
\bar{H}) - (h_{ \lfloor({\varepsilon\log
n})/({4\log|A|}) \rfloor} - \bar{H})
\nonumber\\[-9pt]\\[-9pt]
&\ge&(h_{k_n-1} - \bar{H}) - \frac{1}{\sqrt{n}} \ge\frac{ 3\max(\sqrt
{n},(|A|-1) \operatorname{pen}(n)) }{
n^{1-\varepsilon} },\quad
\nonumber
\end{eqnarray}
where we used that $h_k$ is non-increasing. Comparing (\ref
{eqhpmllb}) to (\ref{eqhpmlub}), the right of (\ref
{eqsplitpml}) is an empty set, and (\ref{eqsplit2}) yields
\[
\Prob\bigl( \hat{k}_{\mathrm{PML}} \bigl(X_1^n\bigr)
< k_n \bigr) \le\Prob\biggl( \overline{B_n \biggl( {
\frac
{\varepsilon
\log n}{4\log|A|} } \biggr)} \biggr) \le12 \RMe ^{1/\RMe } \exp\biggl( -
\frac{7\alpha_0\varepsilon^3}{256\RMe (\alpha
+\alpha_0)} \frac{ n^{\varepsilon/2} }{ \log n } + \frac
{\varepsilon
}{4} \log n \biggr),
\]
if $n\ge(\delta2^{\zeta})^2$, according to Theorem~\ref{thent}.\vadjust{\goodbreak}

(ii) If $\mathrm{IC} = \mathrm{NML}$, by the definition of
the NML information criterion, see Definition~\ref{defNML},
%
%
\begin{eqnarray}
\label{eqsplitnml} && \biggl\{ \mathrm{NML}_{X_1^n} ( m ) \le\mathrm
{NML}_{X_1^n} \biggl( \biggl\lfloor{ \frac{\varepsilon
\log n}{4\log|A|} } \biggr\rfloor
\biggr) \mbox{ for some } m<k_n \biggr\} \cap B_n \biggl(
{ \frac{\varepsilon\log n}{4\log|A|} } \biggr)
\nonumber
\\
&&\quad\subseteq\biggl\{ (n-m)\hat{h}_{m}\bigl(X_1^n
\bigr) - \biggl( n- \biggl\lfloor{ \frac{\varepsilon\log n}{4\log|A|}
} \biggr\rfloor\biggr) \hat
{h}_{ \lfloor({\varepsilon\log n})/({4\log|A|})
\rfloor}\bigl(X_1^n\bigr)
\nonumber
\\
&&\hspace*{16.5pt}\quad\le\log\Sigma\biggl( n, { \biggl\lfloor\frac{\varepsilon\log
n}{4\log|A|} \biggr\rfloor}
\biggr) - \log\Sigma(n,m) \mbox{ for some } m<k_n \biggr\} \cap
B_n \biggl( { \frac{\varepsilon\log n}{4\log|A|} } \biggr)
\nonumber
\\
&&\quad\subseteq\biggl\{ \hat{h}_{m}\bigl(X_1^n
\bigr) - \hat{h}_{ \lfloor
({\varepsilon\log n})/({4\log|A|}) \rfloor}\bigl(X_1^n\bigr)
<
\frac{ \log
\Sigma( n,
{ \lfloor({\varepsilon\log n})/({4\log|A|})
\rfloor} ) }{n- \lfloor(\varepsilon\log n)/(4\log|A|)
\rfloor}
\nonumber\\[-8pt]\\[-8pt]
&&\qquad\hspace*{4pt}\mbox{ for some } m<k_n \biggr\}
\nonumber
\\
&&\qquad{}\cap B_n \biggl( { \frac{\varepsilon\log n}{4\log|A|} } \biggr)
\nonumber
\\
&&\quad\subseteq\biggl\{ h_m - h_{ \lfloor({\varepsilon\log
n})/({4\log
|A|}) \rfloor} <
\frac{ \log\Sigma( n,
{ \lfloor({\varepsilon\log n})/({4\log|A|})
\rfloor} ) }{n- (\varepsilon\log n)/(4\log|A|) } + \frac{2}{n^{
1/2-\varepsilon}} \nonumber\\
&&\qquad\hspace*{4pt}\mbox{ for some } m<k_n \biggr
\},\nonumber
\end{eqnarray}
where in the second relation we used that $\Sigma(n,m)>1$ for any
$m\ge
0$. By Lemma~\ref{lemktml} in the \hyperref[app]{Appendix},
\[
\mathrm{ML}_{k}\bigl(X_1^n\bigr) \le
P_{\mathrm{KT}, k} \bigl(X_1^n\bigr) \exp\biggl(
C_{\mathrm{KT}} |A|^k + \frac{|A|-1}{2} |A|^k \log
\frac{n}{|A|^k} \biggr)
\]
that gives the upper bound
%
%
\begin{equation}
\label{eqsigmaub} \log\Sigma(n,k) \le C_{\mathrm{KT}} |A|^k +
\frac{|A|-1}{2} |A|^k \log\frac{n}{|A|^k}.
\end{equation}
Using (\ref{eqsigmaub}) and (\ref{eqatom}),
\begin{eqnarray*}
&&
\frac{ \log\Sigma( n, { \lfloor
({\varepsilon\log n})/({4\log|A|}) \rfloor} ) }{n-
(\varepsilon\log n)/(4\log|A|) } \\
&&\quad< \biggl( C_{\mathrm{KT}} + \frac
{|A|-1}{2} \log
\frac{n}{|A|^{
\lfloor{ ({\varepsilon\log n})/({4\log|A|}) } \rfloor}} \biggr) \frac
{1}{ n^{1-\varepsilon} }
\\
&&\quad< \biggl( C_{\mathrm{KT}} + \frac{|A|-1}{2} \log n \biggr)
\frac
{1}{ n^{1-\varepsilon} }.
\end{eqnarray*}
Using $\RMe ^x\ge x^2/2+x^4/4!$, $x\ge0$, it follows that $(|A|-1)\log n
\le\sqrt{n}$ if $n\ge24(\log^4 \RMe )(|A|-1)^4$, which implies that
\[
C_{\mathrm{KT}} + \frac{|A|-1}{2} \log n \le\sqrt{n} \qquad\mbox{if }
n\ge
\max\bigl\{ 24\bigl(\log^4 \RMe \bigr) \bigl(|A|-1\bigr)^4,
4C_{\mathrm
{KT}}^2 \bigr\} .
\]
Thus, the expression in (\ref{eqsplitnml}) can be bounded as
%
%
\begin{eqnarray}
\label{eqhnmlub}
&&\frac{ \log\Sigma( n, { \lfloor
({\varepsilon\log n})/({4\log|A|}) \rfloor} ) }{n-
(\varepsilon
\log n)/(4\log|A|) } + \frac{2}{n^{1/2-\varepsilon}}
\nonumber\\[-8pt]\\[-8pt]
&&\quad
< \frac
{3}{n^{1/2-\varepsilon}}
\qquad\mbox{if } n\ge\max\bigl\{ 24\bigl(\log^4 \RMe \bigr)
\bigl(|A|-1\bigr)^4, 4C_{\mathrm
{KT}}^2 \bigr\}.\nonumber
\end{eqnarray}
Now, let $\varepsilon$ and $k_n$ be as in the claim of the proposition.
Then the conditions $h_k-\bar{H} \le\delta2^{-\zeta k}$ and
$\varepsilon\ge(4\log|A|)/\zeta$ imply (\ref{eqdhpmlub}), thus, if
$n\ge(\delta2^{\zeta})^2$, it follows that $k_n \le\frac
{\varepsilon
\log n}{4\log|A|}$, and for any $m<k_n$
%
%
\begin{eqnarray}
\label{eqhnmllb} h_m - h_{ \lfloor({\varepsilon\log n})/({4\log|A|})
\rfloor} &\ge&(h_{k_n-1} -
\bar{H}) - (h_{ \lfloor({\varepsilon\log
n})/({4\log|A|}) \rfloor} - \bar{H})\nonumber\\
&\ge&(h_{k_n-1} - \bar{H}) -
\frac{1}{\sqrt{n}} \\
&\ge&\frac{3}{n^{1/2-\varepsilon}},\nonumber
\end{eqnarray}
where we used that $h_k$ is non-increasing. Comparing (\ref
{eqhnmllb}) to (\ref{eqhnmlub}), the right of (\ref
{eqsplitnml}) is an empty set, and (\ref{eqsplit2}) yields
\[
\Prob\bigl( \hat{k}_{\mathrm{NML}} \bigl(X_1^n\bigr)
< k_n \bigr) \le\Prob\biggl( \overline{B_n \biggl( {
\frac
{\varepsilon
\log n}{4\log|A|} } \biggr)} \biggr)
\le 12 \RMe ^{1/\RMe } \exp\biggl( -
\frac{7\alpha_0\varepsilon^3}{256\RMe (\alpha
+\alpha_0)} \frac{ n^{\varepsilon/2} }{ \log n } + \frac
{\varepsilon
}{4} \log n \biggr),
\]
if $n\ge\max\{ 24(\log^4 \RMe )(|A|-1)^4, 4C_{\mathrm{KT}}^4,
(\delta
2^{\zeta})^2 \}$, according to Theorem~\ref{thent}.

(iii) If $\mathrm{IC} = \mathrm{KT}$, by the definition of
the KT information criterion, see Definition~\ref{defKT}, and using
that $P_{\mathrm{KT}, m} (X_1^n) \le\mathrm{ML}_{m}(X_1^n)$ for any
$0\le m <n$,
%
%
\begin{eqnarray}
\label{eqsplitkt3} && \biggl\{ \mathrm{KT}_{X_1^n} ( m ) \le\mathrm
{KT}_{X_1^n} \biggl( \biggl\lfloor{ \frac{\varepsilon\log
n}{4\log|A|} } \biggr\rfloor
\biggr) \mbox{ for some } m<k_n \biggr\} \cap B_n \biggl(
{ \frac{\varepsilon\log n}{4\log|A|} } \biggr)
\nonumber
\\[-0.8pt]
&&\quad\subseteq\biggl\{ (n-m)\hat{h}_{m}\bigl(X_1^n
\bigr) - \biggl( n- \biggl\lfloor{ \frac{\varepsilon\log n}{4\log|A|}
} \biggr\rfloor\biggr) \hat
{h}_{ \lfloor({\varepsilon\log n})/({4\log|A|})
\rfloor}\bigl(X_1^n\bigr)
\nonumber
\\[-0.8pt]
&&\hspace*{16.4pt}\quad\le\log\mathrm{ML}_{ \lfloor({\varepsilon\log n})/({4\log
|A|}) \rfloor}\bigl(X_1^n
\bigr) - \log P_{\mathrm{KT}, \lfloor
({\varepsilon\log n})/({4\log|A|}) \rfloor} \bigl(X_1^n\bigr)
\mbox{ for some } m<k_n \biggr\} \nonumber\\[-0.8pt]
&&\qquad{}\cap B_n \biggl( {
\frac{\varepsilon\log n}{4\log|A|} } \biggr)
\nonumber
\\[-0.8pt]
&&\quad\subseteq\biggl\{\hat{h}_{m}\bigl(X_1^n
\bigr) - \hat{h}_{ \lfloor({\varepsilon\log n})/({4\log|A|})
\rfloor}\bigl(X_1^n\bigr)
\nonumber\\[-0.8pt]
&&\hspace*{15.6pt}\quad\le\frac{ \log\mathrm{ML}_{ \lfloor({\varepsilon\log
n})/({4\log|A|}) \rfloor}(X_1^n) - \log P_{\mathrm{KT},
\lfloor({\varepsilon\log n})/({4\log|A|}) \rfloor} (X_1^n)
}{n- \lfloor(\varepsilon\log n)/(4\log|A|)
\rfloor} \\[-0.8pt]
&&\qquad\hspace*{4.5pt}\mbox{ for some } m<k_n \biggr\} \nonumber\\[-0.8pt]
&&\qquad{}\cap
B_n \biggl( { \frac{\varepsilon\log n}{4\log|A|} } \biggr)
\nonumber
\\[-0.8pt]
&&\quad\subseteq\biggl\{ h_m - h_{ \lfloor({\varepsilon\log
n})/({4\log
|A|}) \rfloor} \nonumber\\[-0.8pt]
&&\hspace*{15.6pt}\quad\le
\frac{ \log\mathrm{ML}_{ \lfloor
({\varepsilon\log n})/({4\log|A|}) \rfloor}(X_1^n) - \log
P_{\mathrm{KT}, \lfloor({\varepsilon\log n})/({4\log|A|})
\rfloor} (X_1^n) }{n- (\varepsilon\log
n)/(4\log|A|) } + \frac{2}{n^{1/2-\varepsilon}}
\nonumber
\\[-0.8pt]
&&\qquad\hspace*{4.5pt} \mbox{ for some } m<k_n \biggr\}.\nonumber
\end{eqnarray}
By Lemma~\ref{lemktml} in the \hyperref[app]{Appendix},
\begin{eqnarray*}
&& \log\mathrm{ML}_{ \lfloor({\varepsilon\log n})/({4\log|A|})
\rfloor}\bigl(X_1^n\bigr) -
\log P_{\mathrm{KT}, \lfloor
({\varepsilon\log n})/({4\log|A|}) \rfloor} \bigl(X_1^n\bigr)
\\[-0.8pt]
&&\quad\le C_{\mathrm{KT}} |A|^{({\varepsilon\log n})/({4\log|A|})} +
\frac
{|A|-1}{2}
|A|^{({\varepsilon\log n})/({4\log|A|})} \log\frac
{n}{|A|^{ \lfloor({\varepsilon\log n})/({4\log|A|})
\rfloor}},
\end{eqnarray*}
and the proof continues in the same way as in the NML case
(ii).
\end{pf*}

Now, we are ready to prove Theorem~\ref{thpmllogSMP}. We
prove the following theorem that formulates Theorem~\ref{thpmllogSMP}
with explicit constants.
%
%
\begin{theorem}\label{thpmllog}
For any weakly non-null stationary ergodic process with continuity
rates $\bar{\gamma}(k) \le\delta_1 2^{-\zeta_1 k}$ and $\underbar
{\gamma}(k) \ge\delta_2 2^{-\zeta_2 k}$ for some $\zeta_1, \zeta_2,
\delta_1, \delta_2 >0$ ($\zeta_2\ge\zeta_1$), if
\[
\frac{6\log|A|}{\zeta_1} \le\varepsilon< \frac12 ,
\]

\begin{longlist}
\item
the PML Markov order estimator $\hat{k}_{\mathrm{PML}} (X_1^n)$
satisfies that
\[
\Prob\bigl( \hat{k}_{\mathrm{PML}} \bigl(X_1^n\bigr)
\le k_n \bigr) \le12 \RMe ^{1/\RMe } \exp\biggl( -
\frac{7\alpha_0\varepsilon^3}{256\RMe (\alpha+\alpha_0)} \frac{
n^{\varepsilon/2} }{ \log n } + \frac{\varepsilon}{4} \log n \biggr),
\]
if $n\ge(36 \delta_1^{4/3} 2^{(4\zeta_1)/3} \log^2 |A|)/(\log^2
\RMe )$, where
\[
k_n = \frac{1}{2\zeta_2} \biggl( 2\log\delta_2 -3 +
\biggl( \frac12 -\varepsilon\biggr) \log n - \log\max\biggl\{ 1, \bigl(|A|-1\bigr)
\frac
{\operatorname{pen}(n)}{\sqrt{n}} \biggr\} \biggr) ;
\]

\item the Markov order estimator $\hat{k}_{\mathrm{IC}} (X_1^n)$,
where IC is either NML or KT, satisfies that
\[
\Prob\bigl( \hat{k}_{\mathrm{IC}} \bigl(X_1^n\bigr)
\le k_n \bigr) \le12 \RMe ^{1/\RMe } \exp\biggl( -
\frac{7\alpha_0\varepsilon^3}{256\RMe (\alpha+\alpha_0)} \frac{
n^{\varepsilon/2} }{ \log n } + \frac{\varepsilon}{4} \log n \biggr),
\]
if $n\ge\max^2 \{ \sqrt{24} (\log^2 \RMe ) (|A|-1)^2, 2C_{\mathrm{KT}},
(6 \delta_1^{2/3} 2^{(2\zeta_1)/3} \log|A|)/(\log \RMe ) \}$, where
\[
k_n = \frac{1}{2\zeta_2} \biggl( 2\log\delta_2 -3 +
\biggl( \frac12 -\varepsilon\biggr) \log n \biggr) .
\]
[Here, $C_{\mathrm{KT}}$ is the constant in the well-known bound of
$\log\mathrm{ML}_{k}(X_1^n) -\log P_{\mathrm{KT}, k} (X_1^n)$, see
Lemma~\textup{\ref{lemktml}} in the \hyperref[app]{Appendix}.]
\end{longlist}
\end{theorem}
\begin{pf}
By Remark~\ref{remgammaalpha}, $\sum_{k=0}^{+\infty} \bar{\gamma}(k)
\le\sum_{k=0}^{+\infty} \delta_1 2^{-\zeta_1 k} < +\infty$
implies the
$\alpha$-summability. The deviation of the conditional entropies from
the entropy rate will also be controlled by the continuity rates of the
process, and Proposition~\ref{thpmlh} will yield the claim of the theorem.

First, for any $k\le m$,
%
%
\begin{eqnarray}
\label{eqgsplit1}
&& h_k - h_m
\nonumber
\\
&&\quad= \sum_{a\in A} \sum
_{a_{m-k+1}^m \in A^k} \biggl( -P\bigl(a_{m-k+1}^m a\bigr)
\log\frac{ P(a_{m-k+1}^m a) }{
P(a_{m-k+1}^m) } \nonumber\\
&&\qquad\hspace*{60.5pt}{}- \sum_{a_1^{m-k} \in A^{m-k}} -P
\bigl(a_1^m a\bigr) \log\frac{ P(a_1^m a) }{
P(a_1^m) } \biggr)
\\
&&\quad= \sum_{a\in A} \sum
_{a_{m-k+1}^m \in A^k} \biggl( -P\bigl(a_{m-k+1}^m \bigr)
\sum_{a_1^{m-k} \in A^{m-k}} \frac{ P(a_1^m)
}{ P(a_{m-k+1}^m) } \biggl(
\frac{ P(a_{m-k+1}^m a) }{ P(a_{m-k+1}^m) } \log\frac{ P(a_{m-k+1}^m a)
}{ P(a_{m-k+1}^m) } \biggr)
\nonumber
\\
&&\hspace*{58pt}\qquad{} - P\bigl(a_{m-k+1}^m\bigr) \sum
_{a_1^{m-k} \in A^{m-k}} -\frac{
P(a_1^m) }{ P(a_{m-k+1}^m) } \biggl( \frac{ P(a_1^m a) }{ P(a_1^m) }
\log
\frac{ P(a_1^m a) }{ P(a_1^m) } \biggr)\biggr)
\nonumber
\\
&&\quad= \sum_{a_{m-k+1}^m \in A^k} -P\bigl(a_{m-k+1}^m
\bigr) \sum_{a_1^{m-k} \in A^{m-k}} \frac{ P(a_1^m) }{
P(a_{m-k+1}^m) }
\nonumber
\\
&&\qquad\hspace*{0pt}{}\times \sum_{a\in A} \biggl( \frac{ P(a_{m-k+1}^m a) }{
P(a_{m-k+1}^m) } \log
\frac{ P(a_{m-k+1}^m a) }{ P(a_{m-k+1}^m) } - \frac{ P(a_1^m a) }{
P(a_1^m) } \log\frac{ P(a_1^m a) }{ P(a_1^m) } \biggr).\nonumber
\end{eqnarray}
On the right of (\ref{eqgsplit1}), the difference of entropies of the
conditional distributions $\{ P(a| a_{m-k+1}^m), a\in A\}$
and $\{ P(a| a_1^m), a\in A\}$ appears. By Remark~\ref{remgammaeq},
the total variation of these conditional distributions can be upper
bounded as
\[
d_{\TV} \bigl( P\bigl( \cdot| a_{m-k+1}^m\bigr), P
\bigl( \cdot| a_1^m\bigr) \bigr) = \sum
_{a\in A} \bigl\llvert P\bigl(a| a_{m-k+1}^m
\bigr) - P\bigl(a| a_1^m\bigr) \bigr\rrvert\le\bar{
\gamma}(k) .
\]
Hence, applying Lemma~\ref{lemsch} in the \hyperref[app]{Appendix} it
follows, similar
to the bound (\ref{eqsch}) and (\ref{eqsch2}) in the proof of
Proposition~\ref{prent}, that
%
%
\begin{eqnarray}
\label{eqgsplit11} &&\biggl\llvert\sum_{a\in A} P
\bigl(a| a_{m-k+1}^m\bigr) \log P\bigl(a| a_{m-k+1}^m
\bigr) - \sum_{a\in A} P\bigl(a| a_1^m
\bigr) \log P\bigl(a| a_1^m\bigr) \biggr\rrvert
\nonumber
\\
&&\quad\le\frac{\log|A|}{\log \RMe } d_{\TV} \bigl( P\bigl( \cdot|
a_{m-k+1}^m\bigr) , P\bigl( \cdot| a_1^m
\bigr) \bigr) + \frac{1}{\RMe } \frac{1+\nu}{\nu} d_{\TV}^{{1}/({1+\nu})}
\bigl( P\bigl( \cdot| a_{m-k+1}^m\bigr), P\bigl( \cdot|
a_1^m\bigr) \bigr)\qquad
\nonumber\\[-8pt]\\[-8pt]
&&\quad\le\frac{\log|A|}{\log \RMe } \bar{\gamma}(k) + \frac{1}{\RMe }
\frac{1+\nu}{\nu} \bar{\gamma}(k)^{
{1}/({1+\nu})}
\nonumber
\\
&&\quad\le\frac{2\log|A|}{\log \RMe } \frac{1+\nu}{\nu} \bar{\gamma
}(k)^{{1}/({1+\nu})}\nonumber
\end{eqnarray}
for any $\nu>0$, if $\bar{\gamma}(k) \le1/\RMe $. Setting $\nu=1/2$,
combining (\ref{eqgsplit11}) with (\ref{eqgsplit1}) and taking
$m\to
+\infty$ yield the bound
%
%
\begin{equation}
\label{eqklb} h_k - \bar{H} \le\frac{6\log|A|}{\log \RMe } \bar{
\gamma}(k)^{2/3},
\end{equation}
if $\bar{\gamma}(k) \le1/\RMe $. Since $h_k - \bar{H} \le h_k \le\log|A|$,
the bound (\ref{eqklb}) is trivial if $\bar{\gamma}(k) > 1/\RMe $. Hence,
using the assumption $\bar{\gamma}(k) \le\delta_1 2^{-\zeta_1 k}$ of
the theorem,
%
%
\begin{equation}
\label{eqhg} h_k - \bar{H} \le\frac{6\log|A|}{\log \RMe }
\delta_1^{2/3} 2^{-{2\zeta
_1} k/3},
\end{equation}
and the assumption $h_k-\bar{H} \le\delta2^{-\zeta k}$ of Proposition
\ref{thpmlh} is satisfied with
\[
\delta= \frac{6\log|A|}{\log \RMe } \delta_1^{2/3} \quad\mbox{and}
\quad\zeta= \frac{2\zeta_1}{3} .
\]
Thus, the constraint $\varepsilon\ge(4\log|A|)/\zeta$ in Proposition
\ref{thpmlh} becomes $\varepsilon\ge(6\log|A|)/\zeta_1$, and
$n\ge
(\delta2^{\zeta})^2$ becomes
\[
n\ge\frac{36\log^2 |A|}{\log^2 \RMe } \delta_1^{4/3} 2^{(4\zeta_1)/3} .
\]

Next, for any $k < +\infty$,
%
%
\begin{eqnarray}
\label{eqgsplit2} &&h_k - \bar{H}
\nonumber
\\
&&\quad=\sum_{x_{-k}^{-1} \in A^k} \sum
_{a\in A} -P\bigl(x_{-k}^{-1} a\bigr) \log P
\bigl(a| x_{-k}^{-1} \bigr) \nonumber\\
&&\qquad{}+ \int_{A^{\infty}}
\sum_{a\in A} P\bigl(a| x_{-\infty}^{-1}
\bigr) \log P\bigl(a| x_{-\infty}^{-1} \bigr) \mrmd P
\bigl(x_{-\infty}^{-1}\bigr)
\\
&&\quad=\int_{A^{\infty}} \sum_{a\in A} P
\bigl(a| x_{-\infty}^{-1} \bigr) \log\frac{P(a| x_{-\infty}^{-1}
)}{P(a| x_{-k}^{-1} )} \mrmd P
\bigl(x_{-\infty}^{-1}\bigr)
\nonumber
\\
&&\quad=\int_{A^{\infty}} D \bigl( P\bigl( \cdot|
x_{-\infty}^{-1} \bigr) \| P\bigl( \cdot| x_{-k}^{-1}
\bigr) \bigr) \mrmd P\bigl(x_{-\infty}^{-1}\bigr),\nonumber
\end{eqnarray}
where $D( \cdot\| \cdot)$ denotes the Kullback--Leibler
divergence. Using Pinsker's inequality~\cite{Cover,CSbook}, (\ref
{eqgsplit2}) can be lower bounded by
%
%
\begin{equation}
\label{eqgsplit22} \int_{A^{\infty}} \frac12 \biggl( \sum
_{a\in A} \bigl\llvert P\bigl(a| x_{-\infty}^{-1}
\bigr) - P\bigl(a| x_{-k}^{-1} \bigr) \bigr\rrvert
\biggr)^2 \mrmd P\bigl(x_{-\infty}^{-1}\bigr) \ge
\frac12 \underbar{\gamma}(k)^2 \ge\delta_2^2
2^{-2\zeta_2 k -1},
\end{equation}
where in the last inequality we used the assumption $\underbar{\gamma
}(k) \ge\delta_2 2^{-\zeta_2 k}$ of the theorem. Hence, in case (i)
\begin{eqnarray*}
&&\min\biggl\{ k\ge0\dvtx h_k - \bar{H} < \frac{ 4\max(\sqrt
{n},(|A|-1)\operatorname{pen}(n)) }{ n^{1-\varepsilon} }
\biggr\}
\\
&&\quad\ge\min\biggl\{ k\ge0\dvtx\delta_2^2
2^{-2\zeta_2 k-1} < \frac{
4\max
(\sqrt{n},(|A|-1)\operatorname{pen}(n)) }{ n^{1-\varepsilon} } \biggr\}
\\
&&\quad= \min\biggl\{ k\ge0\dvtx k > \frac{1}{2\zeta_2} \bigl( 2\log
\delta_2 -3 + (1-\varepsilon) \log n - \log\max\bigl(\sqrt{n},\bigl(|A|-1\bigr)
\operatorname{pen}(n)\bigr) \bigr) \biggr\}
\\
&&\quad= 1+ \biggl\lfloor\frac{1}{2\zeta_2} \biggl( 2\log\delta_2 -3
+ \biggl( { \frac12} -\varepsilon\biggr) \log n - \log\max\biggl\{ 1, \bigl(|A|-1\bigr)
\frac{\operatorname{pen}(n)}{\sqrt{n}} \biggr\} \biggr) \biggr\rfloor,
\end{eqnarray*}
while in case (ii)
\begin{eqnarray*}
&&\min\biggl\{ k\ge0\dvtx h_k - \bar{H} < \frac{4}{n^{
1/2-\varepsilon}} \biggr
\}
\\
&&\quad\ge\min\biggl\{ k\ge0\dvtx\delta_2^2
2^{-2\zeta_2 k-1} < \frac
{4}{n^{1/2-\varepsilon}} \biggr\}
\\
&&\quad= \min\biggl\{ k\ge0\dvtx k > \frac{1}{2\zeta_2} \biggl( 2\log
\delta_2 -3 + \biggl( { \frac12} -\varepsilon\biggr) \log n \biggr)
\biggr\}
\\
&&\quad= 1+ \biggl\lfloor\frac{1}{2\zeta_2} \biggl( 2\log\delta_2 -3
+ \biggl( { \frac12} -\varepsilon\biggr) \log n \biggr) \biggr\rfloor,
\end{eqnarray*}
and the proof is completed.
\end{pf}

Finally, we prove the following proposition that directly implies
Theorem~\ref{thoracle}.
%
%
\begin{proposition}\label{proporacle}
For any weakly non-null stationary ergodic process with continuity rate
$\bar{\gamma}(k) \le\delta2^{-\zeta k}$, $\zeta, \delta>0$, and for
any $\xi>0$, if $\varepsilon>0$ is so small and $\zeta>0$ is so
large that
\[
\frac12 + \varepsilon< \kappa< 1-\frac{\varepsilon}{4}
\]
and
\[
\frac{6\log|A|}{\zeta} \le\frac{\varepsilon}{1-\kappa} < 2\xi,
\]
the PML Markov order estimator $\hat{k}_{\mathrm{PML}} (X_1^n)$ with
$\operatorname{pen}(n)=n^{\kappa}$ satisfies that
\[
\Prob\biggl( \biggl\llvert\frac{ \hat{k}_{\mathrm{PML}} (X_1^n) }{
k_{\mathrm
{PML},n} } -1 \biggr\rrvert>\xi\biggr)
\le\exp\biggl( -\frac{c_2'\varepsilon^3}{\log n} n^{\varepsilon/2}
\biggr),
\]
if $n$ is sufficiently large, where $c_2'>0$ is a constant depending
only on the distribution of the process.
\end{proposition}
\begin{pf}
The proof of Theorem~\ref{thpmllog} begins with the observation that
the summability of the continuity rate implies the $\alpha
$-summability. Hence, the conditions of Theorem~\ref{thent} are
satisfied now. Moreover, according to (\ref{eqhg}), $\bar{\gamma}(k)
\le\delta2^{-\zeta k}$ also implies that
%
%
\begin{equation}
\label{eqhg2} h_k - \bar{H} \le\frac{6\log|A|}{\log \RMe }
\delta^{2/3} 2^{-{2\zeta}k/3}.
\end{equation}
Set $\xi$, $\varepsilon$ and $\kappa$ as in the conditions of the
proposition, and
define a sequence $k_n\in\mathbb{N}$ such that for sufficiently large
$n$
\begin{eqnarray*}
&&\mbox{\hphantom{ii}(i)\quad} h_{\lfloor(1-\xi/2) k_n \rfloor} - h_{k_n} \ge|A|
n^{-1+\kappa
+\varepsilon/4},\\
&&\mbox{\hphantom{i}(ii)\quad}h_{k_n} - \bar{H} \le\frac{1}{2} \bigl(|A|-1\bigr)^2
n^{-1+\kappa},\\
&&\mbox{(iii)\quad} k_n\le\frac{\varepsilon\log n}{4\log|A|} .
\end{eqnarray*}
Due to (\ref{eqhg2}), such a sequence exists. Since $h_k - \bar{H}$ is
non-negative decreasing, it is sufficient to show this when $h_k - \bar
{H} = \frac{6\log|A|}{\log \RMe } \delta^{2/3} 2^{-{2\zeta}k/{3}}$.
Then, writing $k_n$ in the form $k_n=\nu\log n$,
\[
h_{\lfloor(1-\xi/2) k_n \rfloor} - h_{k_n} \ge\frac12 \frac{6\log
|A|}{\log \RMe }
\delta^{2/3} n^{-
{2\zeta
}(1-\xi/2)\nu/{3}} \quad\mbox{and}\quad h_{k_n} -
\bar{H} = \frac{6\log|A|}{\log \RMe } \delta^{2/3} n^{-{2\zeta}
\nu/{3}},
\]
if $n$ is sufficiently large, that implies (i) and (ii) if
\[
1-\kappa< \frac{2\zeta\nu}{3} < \biggl( 1-\kappa-\frac
{\varepsilon}{4} \biggr)
\frac{1}{1-\xi/2} .
\]
Such $\nu>0$ exists because it follows from the condition $\varepsilon
/(1-\kappa) < 2\xi$ that $1-\kappa< ( 1-\kappa-\varepsilon/4 ) /
(1-\xi
/2)$. Moreover, the condition $\varepsilon/(1-\kappa) \ge(6 \log
|A|)/\zeta$ implies $\nu\le\varepsilon/(4 \log|A|)$ satisfying (iii).

First, recall the definition of $B_n ( \frac{\varepsilon\log
n}{4\log|A|} )$ in (\ref{eqBndef}) in the proof of Proposition
\ref{thpmlh}. Similar to (\ref{eqsplit2}) and (\ref{eqsplitpml}), we
can write that
%
%
\begin{eqnarray}
\label{eqolowerb} && \bigl\{ \hat{k}_{\mathrm{PML}} \bigl(X_1^n
\bigr) < (1-\xi/2) k_n \bigr\} \cap B_n \biggl( {
\frac{\varepsilon\log n}{4\log|A|} } \biggr)
\\
&&\quad\subseteq\bigl\{ \mathrm{PML}_{X_1^n} ( m ) \le
\mathrm{PML}_{X_1^n} ( k_n ) \mbox{ for some } m<(1-\xi/2)
k_n \bigr\} \cap B_n \biggl( { \frac{\varepsilon\log n}{4\log|A|} }
\biggr)
\nonumber
\\
&&\quad= \bigl\{ \mathrm{PML}_{o,n} ( m ) + (n-m) \bigl(\hat
{h}_{m}\bigl(X_1^n\bigr) - h_m
\bigr) \le\mathrm{PML}_{o,n} ( k_n ) + (n-k_n)
\bigl(\hat{h}_{k_n}\bigl(X_1^n\bigr) -
h_{k_n}\bigr)
\nonumber
\\
&&\hspace*{25pt} \mbox{ for some } m<(1-\xi/2) k_n \bigr\} \cap B_n
\biggl( { \frac{\varepsilon\log n}{4\log|A|} } \biggr)
\nonumber
\\
\label{eqolowerhalf} &&\quad\subseteq\biggl\{ \mathrm{PML}_{o,n} ( m
) - \frac
{2n}{n^{1/2-\varepsilon}} \le\mathrm{PML}_{o,n} ( k_n )
\mbox{ for some } m<(1-\xi/2) k_n \biggr\}
\\
&&\quad\subseteq\biggl\{ (n-m) h_{m} - ( n- k_n )
h_{k_n} \le\bigl(|A|-1\bigr) \bigl( |A|^{k_n} - |A|^m \bigr)
\operatorname{pen}(n) + \frac
{2n}{n^{1/2-\varepsilon}}
\nonumber
\\
&&\hspace*{25pt} \mbox{ for some } m<(1-\xi/2) k_n \biggr\}
\nonumber
\\
&&\quad\subseteq\biggl\{ h_m - h_{ k_n } \le
\frac{(|A|-1) |A|^{
({\varepsilon\log n})/({4\log|A|})} \operatorname{pen}(n)}{n-
(\varepsilon\log n)/(4\log
|A|) } + \frac{2}{n^{1/2-\varepsilon}} \mbox{ for some } m<(1-\xi/2)
k_n \biggr\}
\nonumber
\\
\label{eqolower1} &&\quad\subseteq\bigl\{ h_{\lfloor(1-\xi/2) k_n
\rfloor} - h_{ k_n }
< |A| n^{-1+\kappa+\varepsilon/4} \bigr\}
\end{eqnarray}
that is empty set by (i), if $n$ is large enough and $k_n \le\frac
{\varepsilon\log n}{4\log|A|}$. The latter is satisfied because
of~\textup{(iii)}. On the other hand,
%
%
\begin{eqnarray}
\label{eqoupperb} && \bigl\{ \hat{k}_{\mathrm{PML}} \bigl(X_1^n
\bigr) > (1+\xi/2) k_n \bigr\} \cap B_n \biggl( {
\frac{\varepsilon\log n}{4\log|A|} } \biggr) \cap\biggl\{ \hat
{k}_{\mathrm{PML}}
\bigl(X_1^n\bigr) \le\frac{\varepsilon
\log
n}{4\log|A|} \biggr\}
\\[-0.5pt]
&&\quad\subseteq\biggl\{ \mathrm{PML}_{X_1^n} ( m ) < \mathrm
{PML}_{X_1^n} ( k_n ) \mbox{ for some } (1+\xi/2)
k_n <m\le\frac{\varepsilon\log n}{4\log|A|} \biggr\}\nonumber\\[-0.5pt]
&&\qquad{} \cap B_n \biggl( {
\frac{\varepsilon\log n}{4\log|A|} } \biggr)
\nonumber
\\[-0.5pt]
\label{eqoupperhalf} &&\quad\subseteq\biggl\{ \mathrm{PML}_{o,n} ( m
) - \frac
{2n}{n^{1/2-\varepsilon}} < \mathrm{PML}_{o,n} ( k_n ) \mbox{
for some } m>(1+\xi/2) k_n \biggr\}
\\[-0.5pt]
&&\quad\subseteq\biggl\{ \bigl(|A|-1\bigr) \bigl( |A|^m - |A|^{k_n}
\bigr) \operatorname{pen}(n) - \frac{2n}{n^{1/2-\varepsilon}} < (
n- k_n )
h_{k_n} - (n-m) h_{m}
\nonumber
\\[-0.5pt]
&&\hspace*{25pt} \mbox{ for some } m>(1+\xi/2) k_n \biggr\}
\nonumber
\\[-0.5pt]
&&\quad\subseteq\biggl\{ \bigl(|A|-1\bigr) \bigl( |A|^m - |A|^{k_n}
\bigr) \frac
{\operatorname{pen}(n)}{n} - \frac{2}{n^{1/2-\varepsilon}} <
h_{k_n} - \biggl(1-
\frac{m}{n} \biggr) h_{m}
\nonumber
\\[-0.5pt]
&&\hspace*{25pt} \mbox{ for some } m>(1+\xi/2) k_n \biggr\}
\nonumber
\\[-0.5pt]
&&\quad\subseteq\biggl\{ \bigl(|A|-1\bigr) \bigl( |A|^m - |A|^{k_n}
\bigr) \frac
{\operatorname{pen}(n)}{n} - \frac{2}{n^{1/2-\varepsilon}} - \frac
{m}{n}\bar{H} \nonumber\\[-0.5pt]
&&\hspace*{26.5pt}
< (
h_{k_n} - \bar{H} ) - \biggl(1-\frac{m}{n} \biggr) (
h_{k_n} - \bar{H} )
\nonumber
\\[-0.5pt]
&&\hspace*{25pt} \mbox{ for some } m>(1+\xi/2) k_n \biggr\}
\nonumber
\\[-0.5pt]
\label{eqoupper1} &&\quad\subseteq\biggl\{ h_{k_n} - \bar{H} >
\frac{ (|A|-1)^2 }{2} n^{-1+\kappa} \biggr\}
\end{eqnarray}
that is empty set by (ii), if $n$ is large enough.

Observe that
%
%
\begin{equation}
\label{eqknko} \frac{1+\xi/2}{1+\xi} k_n \le k_{\mathrm{PML},n} \le
\frac{1-\xi
/2}{1-\xi} k_n,
\end{equation}
if $n$ is sufficiently large. Indeed, on indirect way the following
sequence of implications can be written
\begin{eqnarray*}
k_{\mathrm{PML},n} < \frac{1+\xi/2}{1+\xi} k_n \quad&\Rightarrow\quad&
k_{\mathrm{PML},n} < 1-\frac{\xi/2}{1+\xi} k_n \quad\Rightarrow\quad
k_{\mathrm{PML},n} < (1-\xi/2) k_n
\\
&\Rightarrow\quad&\mathrm{PML}_{o,n} ( m ) < \mathrm{PML}_{o,n}
( k_n ) \qquad\mbox{for some } m<(1-\xi/2) k_n
\\
&\Rightarrow\quad&\mathrm{PML}_{o,n} ( m ) - \frac
{2n}{n^{1/2-\varepsilon}} <
\mathrm{PML}_{o,n} ( k_n ) \\
&&\mbox{for some } m<(1-
\xi/2) k_n
\end{eqnarray*}
that does not hold by (\ref{eqolowerhalf}) and (\ref{eqolower1}) if
$n$ is large enough, and
\begin{eqnarray*}
k_{\mathrm{PML},n} > \frac{1-\xi/2}{1-\xi} k_n \quad&\Rightarrow&\quad
k_{\mathrm{PML},n} > 1+\frac{\xi/2}{1-\xi} k_n \Rightarrow
k_{\mathrm{PML},n} > (1+\xi/2) k_n
\\
&\Rightarrow&\quad\mathrm{PML}_{o,n} ( m ) < \mathrm{PML}_{o,n}
( k_n ) \qquad\mbox{for some } m>(1+\xi/2) k_n
\\
&\Rightarrow&\quad\mathrm{PML}_{o,n} ( m ) - \frac
{2n}{n^{1/2-\varepsilon}} <
\mathrm{PML}_{o,n} ( k_n ) \\
&&\quad\mbox{for some } m>(1+
\xi/2) k_n
\end{eqnarray*}
that does not hold either by (\ref{eqoupperhalf}) and (\ref
{eqoupper1}) if $n$ is large enough.

Finally, using (\ref{eqknko}), we get
\begin{eqnarray*}
&& \Prob\biggl( \biggl\llvert\frac{ \hat{k}_{\mathrm{PML}} (X_1^n) }{
k_{\mathrm{PML},n} } -1 \biggr\rrvert>\xi\biggr)
\\
&&\quad=\Prob\bigl( \hat{k}_{\mathrm{PML}} \bigl(X_1^n
\bigr) < (1-\xi) k_{\mathrm
{PML},n} \bigr) + \Prob\bigl( \hat{k}_{\mathrm{PML}}
\bigl(X_1^n\bigr) > (1+\xi) k_{\mathrm
{PML},n} \bigr)
\\
&&\quad\le\Prob\bigl( \hat{k}_{\mathrm{PML}} \bigl(X_1^n
\bigr) < (1-\xi/2) k_n \bigr) + \Prob\bigl( \hat{k}_{\mathrm{PML}}
\bigl(X_1^n\bigr) > (1+\xi/2) k_n \bigr)
\\
&&\quad\le\Prob\biggl( \bigl\{ \hat{k}_{\mathrm{PML}} \bigl(X_1^n
\bigr) < (1-\xi/2) k_n \bigr\} \cap B_n \biggl( {
\frac{\varepsilon\log
n}{4\log
|A|} } \biggr) \biggr)
\\
&&\qquad{} + \Prob\biggl( \bigl\{ \hat{k}_{\mathrm{PML}} \bigl(X_1^n
\bigr) > (1+\xi/2) k_n \bigr\} \cap B_n \biggl( {
\frac{\varepsilon\log
n}{4\log
|A|} } \biggr) \cap\biggl\{ \hat{k}_{\mathrm{PML}}
\bigl(X_1^n\bigr) \le\frac
{\varepsilon\log n}{4\log|A|} \biggr\}
\biggr)
\\
&&\qquad{} +2\Prob\biggl( \overline{ B_n \biggl( {
\frac
{\varepsilon\log n}{4\log|A|} } \biggr) } \biggr) + \Prob\biggl( \hat
{k}_{\mathrm{PML}}
\bigl(X_1^n\bigr) > \frac{\varepsilon
\log
n}{4\log|A|} \biggr),
\end{eqnarray*}
where the first two terms are zero if $n$ is large enough by (\ref
{eqolowerb})--(\ref{eqolower1}) and (\ref{eqoupperb})--(\ref
{eqoupper1}). Using Proposition~\ref{propqminupper} with $r_n=n-1$,
$k_n= \lfloor\frac{\varepsilon\log n}{4\log|A|}
\rfloor$
and $m_n= \lfloor\frac{\varepsilon\log n}{6\log|A|}
\rfloor$,
\[
\Prob\biggl( \hat{k}_{\mathrm{PML}} \bigl(X_1^n
\bigr) > \frac{\varepsilon
\log
n}{4\log|A|} \biggr) \le\exp\bigl( - \mathcal{O}
\bigl(n^{\kappa+\varepsilon/4}\bigr) \bigr),
\]
because
\[
n \bar{\gamma}(m_n) \le n \delta2^{\zeta} \exp\biggl( -\zeta
\frac{\varepsilon\log
n}{6\log|A|} \biggr) = \delta2^{\zeta} n^{1-\zeta\varepsilon/(6\log|A|)},
\]
but $1-\zeta\varepsilon/(6\log|A|) < \kappa+\varepsilon/4$
according to
the condition $\varepsilon/(1-\kappa) \ge(6 \log|A|)/\zeta$. Then the
claim of the proposition follows from Theorem~\ref{thent}.
\end{pf}

\section{Process estimation proofs}\label{secproofappl}

In this section, we consider the estimation of stationary ergodic
processes by finite memory processes. First, define
\[
\beta_1 = \frac{1}{ \prod_{j=1}^{+\infty} (1-2\bar{\gamma}(j)) }
\]
and
\[
\beta_2 = \sup_{k\ge1} 2|A|\frac{1- (1-2|A|\bar{\gamma}(k))^k}{
k \bar
{\gamma}(k) \prod_{j=1}^{+\infty} (1-2|A|\bar{\gamma}(j))^2 } .
\]
Clearly, if $\sum_{k=1}^{\infty} \bar{\gamma}(k)<+\infty$, then
$\beta_1,\beta_2<+\infty$.

Now we prove the following theorem that formulates Theorem \ref
{thqminSMP} with explicit constants.
%
%
\begin{theorem}\label{thqmin}
For any non-null stationary ergodic process with summable continuity
rate and uniformly convergent restricted continuity rate with
parameters $\theta_1$, $\theta_2$, $k_{\theta}$, for any $\mu_n>0$,
the empirical Markov estimator of the process with the order estimated
by the bounded PML Markov order estimator $\hat{k}_{\mathrm{PML}}
(X_1^n | \eta\log n)$, $\eta>0$, with penalty function
$\operatorname{pen}(n)\le
\mathcal{O}(\sqrt{n})$ satisfies
\begin{eqnarray*}
&& \Prob\biggl( \bar{d} \bigl( X_1^n, \hat{X} \bigl[
\hat{k}_{\mathrm
{PML}} \bigl(X_1^n | \eta\log n\bigr)
\bigr]_1^n \bigr) > \frac{\beta_2}{\pinf^{2}} g_n +
\frac{1}{n^{1/2-\mu_n}} \biggr)
\\
&&\quad\le2 \RMe ^{1/\RMe } |A|^{K_n+h_n+2} \exp\biggl\{ -
\frac{ \pinf^2 }{ 16 \RMe  |A|^3 (\alpha+ \pinf) (\beta_1+1)^2 } \frac{
(n-K_n-h_n)}{ (1+K_n+h_n) n }
\\
&&\hspace*{87.4pt}\qquad{}\times 4^{-(K_n+h_n) |\log\pinf|} \biggl[ 4^{\mu_n\log n} - \frac
{(K_n+h_n)|\log\pinf|(\beta_1+1)^2}{2} \biggr]
\biggr\}
\\
&&\qquad{}+ 12 \RMe ^{1/\RMe } \exp\biggl( -\frac{7\alpha_0 (\log|A|)^3
\eta^3}{4\RMe (\alpha+\alpha_0)}
\frac{ n^{\eta2\log|A|} }{ \log n } + \bigl(\eta\log|A|\bigr)\log n \biggr)
\\
&&\qquad{}+ \exp\biggl( -\bigl(|A|-1\bigr) |A|^{K_n+h_n+1} \\
&&\hspace*{32pt}\qquad{}\times\operatorname{pen}(n)
\biggl[ 1 - \frac
{1}{|A|^{1+h_n}} - \frac{1}{2\operatorname{pen}(n)} \bigl( \log n -
(K_n+h_n) \log|A| \bigr) \biggr]
\\
&&\hspace*{32pt}\qquad{} + \frac{c \operatorname{pen}(n)}{\pinf/\log \RMe } +
|A|^{K_n+h_n+1} C_{\mathrm{KT}} + \log(
\eta\log n) \biggr),
\end{eqnarray*}
if $n$ is so large that
%
%
\begin{equation}
\label{eqnbound} \min\biggl\{ \biggl\lfloor\frac{\eta}{\theta_2}\log n
\biggr
\rfloor, k\ge0\dvtx\bar{\gamma}(k) < \biggl( \frac{ 6\max( \sqrt{n},(|A|-1)
\operatorname{pen}(n))
}{ \pinf n^{1-\eta\log(|A|^4/\pinf)} }
\biggr)^{1/(2\theta_1)} \biggr\} \ge k_{\theta},
\end{equation}
where
\begin{eqnarray*}
g_n &=& \max\biggl\{ \bar{\gamma} \biggl( \biggl\lfloor
\frac{\eta
}{\theta_2}\log n \biggr\rfloor\biggr), \biggl( \frac{ 6\max( 1,
(|A|-1)
({\operatorname{pen}(n)})/{\sqrt{n}}) }{ \pinf n^{1/2-\eta
\log(|A|^4/\pinf)} }
\biggr)^{1/(2\theta_1)} \biggr\} ,
\\
K_n &=& K_n \biggl( r_n, \bar{\gamma},
\frac{c}{n}\operatorname{pen}(n) \biggr) ,
\end{eqnarray*}
and $c>0$ is an arbitrary constant and $h_n\in\mathbb{N}$ is an
arbitrary sequence.
\end{theorem}

The proof is based on the following two propositions.
%
%
\begin{proposition}\label{propqmin}
For any non-null and $\alpha$-summable stationary ergodic process with
uniformly convergent restricted continuity rate with parameters $\theta
_1$, $\theta_2$, $k_{\theta}$,

\begin{longlist}
\item
the bounded PML Markov order estimator $\hat
{k}_{\mathrm
{PML}} (X_1^n | \eta\log n)$ with penalty function $\operatorname
{pen}(n)\le\mathcal
{O}(\sqrt{n})$ satisfies that
\[
\Prob\bigl( \hat{k}_{\mathrm{PML}} \bigl(X_1^n |
\eta\log n\bigr) < k_n \bigr) \le12 \RMe ^{1/\RMe } \exp\biggl( -
\frac{7\alpha_0 (\log|A|)^3 \eta^3}{4\RMe (\alpha+\alpha_0)} \frac{ n^{\eta
2\log|A|} }{ \log n } + \bigl(\eta\log|A|\bigr)\log n \biggr),
\]
if $n$ is so large that $k_n\ge k_{\theta}$, where
\[
k_n = \min\biggl\{ \biggl\lfloor\frac{\eta}{\theta_2}\log n \biggr
\rfloor, k\ge0\dvtx\bar{\gamma}(k) < \biggl( \frac{ 6\max( \sqrt{n},(|A|-1)
\operatorname{pen}(n))
}{ \pinf n^{1-\eta\log(|A|^4/\pinf)} }
\biggr)^{1/(2\theta_1)} \biggr\} ;
\]

\item the bonded Markov order estimator $\hat{k}_{\mathrm{IC}} (X_1^n
| \eta\log n)$, where IC is either NML or KT, satisfies that
\[
\Prob\bigl( \hat{k}_{\mathrm{IC}} \bigl(X_1^n |
\eta\log n\bigr) < k_n \bigr) \le12 \RMe ^{1/\RMe } \exp\biggl( -
\frac{7\alpha_0 (\log|A|)^3 \eta^3}{4\RMe (\alpha+\alpha_0)} \frac{ n^{\eta
2\log|A|} }{ \log n } + \bigl(\eta\log|A|\bigr)\log n \biggr),
\]
if $n$ is so large that $k_n\ge k_{\theta}$ and $n\ge\max^2 \{
\sqrt{24}(\log^2 \RMe )(|A|-1)^2, 2C_{\mathrm{KT}} \}$, where
\[
k_n = \min\biggl\{ \biggl\lfloor\frac{\eta}{\theta_2}\log n \biggr
\rfloor, k\ge0\dvtx\bar{\gamma}(k) < \biggl( \frac{ 6 }{ \pinf n^{
1/2-\eta\log
(|A|^4/\pinf)} }
\biggr)^{1/(2\theta_1)} \biggr\} .
\]
\end{longlist}
\end{proposition}
\begin{pf}
First, define $B_n( \eta\log n )$ similar to (\ref{eqBndef}) in the
proof of Proposition~\ref{thpmlh}. Similar to (\ref
{eqsplit2})--(\ref
{eqhpmlub}), we can write for any $k_n\le(\eta/\theta_2)\log n$ that
%
%
\begin{eqnarray}
\label{eqsplit22} &&\Prob\bigl( \hat{k}_{\mathrm{PML}} \bigl(X_1^n
| \eta\log n\bigr) < k_n \bigr)
\nonumber
\\
&&\quad\le\Prob\biggl( h_m - h_{\lfloor\eta\log n \rfloor} <
\frac{
3\max
( \sqrt{n},(|A|-1) \operatorname{pen}(n)) }{ n^{1-\eta4 \log
|A|} } \mbox{ for some } m<k_n \biggr)
\\
&&\qquad{} +\Prob\bigl( \overline{B_n( \eta\log n )} \bigr).\nonumber
\end{eqnarray}

Now, the difference $h_m - h_{\lfloor\eta\log n \rfloor}$ in (\ref
{eqsplit22}) is controlled as follows. For any $m\le k$,
%
%
\begin{eqnarray}
\label{eqgsplit22a} &&h_m - h_k
\nonumber
\\
&&\quad=\sum_{a\in A} \sum
_{a_1^k \in A^k} \bigl( -P\bigl(a_1^k a\bigr)
\log P\bigl(a| a_{k-m+1}^k \bigr) + P\bigl(a_1^k
a\bigr) \log P\bigl(a| a_1^k \bigr) \bigr)
\nonumber\\[-8pt]\\[-8pt]
&&\quad=\sum_{a_1^k \in A^k} P\bigl(a_1^k
\bigr) \sum_{a\in A} P\bigl(a| a_1^k
\bigr) \log\frac{P(a| a_1^k )}{P(a|
a_{k-m+1}^k )}
\nonumber
\\
&&\quad=\sum_{a_1^k \in A^k} P\bigl(a_1^k
\bigr) D \bigl( P\bigl( \cdot| a_1^k \bigr) \| P\bigl(
\cdot| a_{k-m+1}^k \bigr) \bigr).\nonumber
\end{eqnarray}
Using Pinsker's inequality~\cite{Cover,CSbook}, (\ref{eqgsplit22a})
can be lower bounded by
%
%
\begin{eqnarray}
\label{eqgsplit222} &&\sum_{a_1^k \in A^k} P
\bigl(a_1^k\bigr) \frac12 \biggl( \sum
_{a\in A} \bigl\llvert P\bigl(a| a_1^k
\bigr) - P\bigl(a| a_{k-m+1}^k \bigr) \bigr\rrvert
\biggr)^2
\nonumber
\\
&&\quad\ge\frac12 \bar{\gamma}(m|k)^2 \min_{a_1^k\in A^k} P
\bigl(a_1^k\bigr)
\\
&&\quad\ge\frac12 \bar{\gamma}(m|k)^2 \pinf^{ k}.\nonumber
\end{eqnarray}
Using (\ref{eqgsplit222}) and the assumption $\bar{\gamma
}(k)^{\theta
_1} \le\bar{\gamma}( k | \lceil\theta_2 k \rceil)$ if
$k\ge k_{\theta}$ ($\theta_1\ge1,\theta_2>1$), it follows that
\[
h_{k} - h_{ \lceil\theta_2 k \rceil} \ge\tfrac12 \bar{\gamma}\bigl( k |
\lceil
\theta_2 k \rceil\bigr)^2 \pinf^{ \lceil\theta_2 k \rceil} \ge
\tfrac12 \bar{\gamma}(k)^{2\theta_1} \pinf^{ \theta_2 k +1} \qquad\mbox
{if } k
\ge k_{\theta} .
\]
Hence, we can write
%
%
\begin{eqnarray}
\label{eqknpml2} &&\min\biggl\{ k\ge k_{\theta}\dvtx h_{k} -
h_{ \lceil\theta_2 k
\rceil} < \frac{ 3\max( \sqrt{n},(|A|-1) \operatorname{pen}(n))
}{ n^{1-\eta4
\log|A|} } \biggr\}
\nonumber
\\
&&\quad\ge\min\biggl\{ k\ge k_{\theta}\dvtx\bar{\gamma}(k) < \biggl(
\frac{
6\max( \sqrt{n},(|A|-1) \operatorname{pen}(n)) }{ n^{1-\eta4
\log|A|} } 2^{-(\theta_2 k +1)\log\pinf} \biggr)^{1/(2\theta_1)} \biggr
\}\quad
\\
&&\quad\ge\min\biggl\{ k\ge k_{\theta}\dvtx\bar{\gamma}(k) < \biggl(
\frac{
6\max( \sqrt{n},(|A|-1) \operatorname{pen}(n)) }{ \pinf
n^{1-\eta\log
(|A|^4/\pinf)} } \biggr)^{1/(2\theta_1)} \biggr\}.\nonumber
\end{eqnarray}
Let $k_n$ be as in the claim of the proposition and suppose that $k_n
\ge k_{\theta}$. Then, since $h_k$ is non-increasing, for any $m<k_n\le
(\eta/\theta_2)\log n$
%
%
\begin{equation}
\label{eqhpmllb2} h_m - h_{\lfloor\eta\log n \rfloor} \ge h_{k_n-1} -
h_{ \lceil\theta_2 (k_n-1) \rceil} \ge\frac{ 3\max( \sqrt{n},(|A|-1)
\operatorname{pen}(n)) }{
n^{1-\eta4 \log
|A|} }.
\end{equation}

Applying (\ref{eqhpmllb2}) to (\ref{eqsplit22}), the first term
on the right in (\ref{eqsplit22}) equals zero, therefore
\begin{eqnarray*}
\Prob\bigl( \hat{k}_{\mathrm{PML}} \bigl(X_1^n\bigr)
< k_n \bigr) &\le&\Prob\bigl( \overline{B_n( \eta\log n )}
\bigr)
\\
&\le&12 \RMe ^{1/\RMe } \exp\biggl( -\frac{7\alpha_0 (\log|A|)^3 \eta
^3}{4\RMe (\alpha+\alpha_0)} \frac{ n^{\eta2\log|A|} }{ \log n } + \bigl(
\eta\log|A|\bigr)\log n \biggr)
\end{eqnarray*}
by Theorem~\ref{thent} with $\varepsilon=\eta4\log|A|$.

In cases $\mathrm{IC}=\mathrm{NML}$ and $\mathrm{IC}=\mathrm{KT}$, the proofs deviate from the above similar to
as (ii) and (iii) deviate from (i) in the proof of Proposition
\ref{thpmlh}. Now, instead of (\ref{eqhnmlub}) we have
\begin{eqnarray*}
&&\frac{ \log\Sigma( n, \lfloor\eta\log n
\rfloor
) }{n- \eta\log n } + \frac{2}{n^{1/2-\eta4\log|A|}}\\
&&\quad< \frac
{3}{n^{1/2-\eta4\log|A|}} \qquad\mbox{if }
n\ge\max\bigl\{ 24\bigl(\log^4 \RMe \bigr) \bigl(|A|-1\bigr)^4,
4C_{\mathrm
{KT}}^2 \bigr\} .
\end{eqnarray*}
\upqed\end{pf}
%
%
\begin{proposition}\label{propqminupper}
For any non-null stationary ergodic process, the bounded PML Markov
order estimator $\hat{k}_{\mathrm{PML}} (X_1^n | r_n)$ satisfies that
\begin{eqnarray*}
&&\Prob\bigl( \hat{k}_{\mathrm{PML}} \bigl(X_1^n |
r_n\bigr) > k_n \bigr)
\\
&&\quad\le\exp\biggl( \log(r_n-k_n) +
\frac{(n-m_n)\bar{\gamma
}(m_n)}{\pinf
/\log \RMe } + \bigl(|A|-1\bigr)|A|^{m_n} \operatorname{pen}(n)
\\
&&\hspace*{19pt}\qquad{} + |A|^{k_n+1} \biggl[ C_{\mathrm{KT}} + \frac{|A|-1}{2}
\log\frac
{n}{|A|^{k_n+1}} - \bigl(|A|-1\bigr) \operatorname{pen}(n) \biggr] \biggr)
\end{eqnarray*}
for any $0\le m_n \le k_n \le r_n \le n$.
\end{proposition}
\begin{pf}
For any $m\ge0$,
%
%
\begin{equation}
\label{eqpfact1} P \bigl(x_1^n\bigr) = P
\bigl(x_1^m\bigr) \prod_{i=m+1}^n
P\bigl(x_i | x_1^{i-1}\bigr) \le\Biggl( \prod
_{i=m+1}^n P\bigl(x_i |
x_{i-m}^{i-1}\bigr) \Biggr) \prod
_{i=m+1}^n \frac{ P(x_i | x_1^{i-1}) }{ P(x_i | x_{i-m}^{i-1}) }.
\end{equation}
Using $P(x_i | x_1^{i-1}) \le P(x_i | x_{i-m}^{i-1}) + \bar{\gamma}(m)$
and $P(x_i | x_{i-m}^{i-1}) \ge\pinf$, (\ref{eqpfact1}) can be upper
bounded by
%
%
\begin{equation}
\label{eqpfact2} \Biggl( \prod_{i=m+1}^n P
\bigl(x_i | x_{i-m}^{i-1}\bigr) \Biggr) \biggl( 1+
\frac{\bar{\gamma}(m)}{\pinf} \biggr)^{n-m} \le\mathrm{ML}_{m}
\bigl(x_1^n\bigr) \biggl( 1+ \frac{\bar{\gamma}(m)}{\pinf}
\biggr)^{n-m}.
\end{equation}
Now, let $C_{n,k}= \{ \hat{k}_{\mathrm{PML}} (X_1^n | r_n) = k
\}$. By the definition of the PML information criterion, see
Definition~\ref{defPML}, for any $0\le m_n,k \le r_n$
%
%
\begin{eqnarray}
\label{eqsh1} \log\mathrm{ML}_{m_n}\bigl(X_1^n
\bigr) &\le&\log\mathrm{ML}_{k}\bigl(X_1^n
\bigr)
- \bigl(|A|-1\bigr)|A|^k \operatorname{pen}(n)\nonumber\\[-8pt]\\[-8pt]
&&{} + \bigl(|A|-1\bigr)|A|^{m_n}
\operatorname{pen}(n) \qquad\mbox{if } X_1^n \in
C_{n,k} .
\nonumber
\end{eqnarray}
By Lemma~\ref{lemktml} in the \hyperref[app]{Appendix},
%
%
\begin{equation}
\label{eqktml} \mathrm{ML}_{k}\bigl(X_1^n
\bigr) \le P_{\mathrm{KT}, k} \bigl(X_1^n\bigr) \exp
\biggl( C_{\mathrm{KT}} |A|^k + \frac{|A|-1}{2} |A|^k
\log\frac{n}{|A|^k} \biggr).
\end{equation}
Combining (\ref{eqpfact2}), (\ref{eqsh1}) and (\ref{eqktml}),
\begin{eqnarray*}
P \bigl(X_1^n\bigr) &\le& P_{\mathrm{KT}, k}
\bigl(X_1^n\bigr) \biggl( 1+ \frac{\bar
{\gamma}(m_n)}{\pinf}
\biggr)^{n-m_n} \\
&&{}\times\exp\biggl( C_{\mathrm{KT}} |A|^k +
\frac{|A|-1}{2} |A|^k \log\frac
{n}{|A|^k}
\\
&&\hspace*{33pt}{}- \bigl(|A|-1\bigr)|A|^k \operatorname{pen}(n) + \bigl(|A|-1\bigr)|A|^{m_n}
\operatorname{pen}(n) \biggr) \qquad\mbox{if } X_1^n
\in C_{n,k} ,
\end{eqnarray*}
that implies
%
%
\begin{eqnarray}
\label{eqsh2} P (C_{n,k}) &\le& \biggl( 1+ \frac{\bar{\gamma
}(m_n)}{\pinf}
\biggr)^{n-m_n} \nonumber\\
&&{}\times\exp\biggl( C_{\mathrm{KT}} |A|^k +
\frac{|A|-1}{2} |A|^k \log\frac
{n}{|A|^k}
\nonumber
\\
&&\hspace*{32.5pt}{} - \bigl(|A|-1\bigr)|A|^k \operatorname{pen}(n) + \bigl(|A|-1\bigr)|A|^{m_n}
\operatorname{pen}(n) \biggr)
\\
&\le&\exp\biggl( \frac{(n-m_n)\bar{\gamma}(m_n)}{\pinf/\log \RMe } +
\bigl(|A|-1\bigr)|A|^{m_n}
\operatorname{pen}(n)
\nonumber
\\
&&\hspace*{21pt}{} + |A|^k \biggl[ C_{\mathrm{KT}} + \frac{|A|-1}{2} \log
\frac{n}{|A|^k} - \bigl(|A|-1\bigr) \operatorname{pen}(n) \biggr]
\biggr),\nonumber
\end{eqnarray}
where in the last inequality we used $\log(1+x)\le x\log \RMe $, $x\ge0$.
In the exponent of (\ref{eqsh2}), it may be assumed that $|A|^k$ is
multiplied by a negative number otherwise the bound is trivial. Then,
the claim of the lemma follows from (\ref{eqsh2}) as
\[
\Prob\bigl( \hat{k}_{\mathrm{PML}} \bigl(X_1^n |
r_n\bigr) > k_n \bigr) \le\sum
_{k=k_n+1}^{r_n} P (C_{n,k})
\le(r_n-k_n) P (C_{n,k_n+1}) .
\]
\upqed\end{pf}

Now, we are ready to prove Theorem~\ref{thqmin}.
\begin{pf*}{Proof of Theorem~\ref{thqmin}}
Letting
\[
G_n = \bigl\{ \bar{\gamma} \bigl( \hat{k}_{\mathrm{PML}}
\bigl(X_1^n | \eta\log n\bigr) \bigr) \le
g_n \bigr\}
\]
and
\[
H_n = \bigl\{ \hat{k}_{\mathrm{PML}} \bigl(X_1^n
| \eta\log n\bigr) \le k_n \bigr\} ,
\]
write
%
%
\begin{eqnarray}
\label{eqfinal} &&\Prob\biggl( \bar{d} \bigl( X_1^n,
\hat{X} \bigl[ \hat{k}_{\mathrm
{PML}} \bigl(X_1^n |
\eta\log n\bigr) \bigr]_1^n \bigr) > \frac{\beta_2}{\pinf^{2}}
g_n + \frac{1}{n^{1/2-\mu_n}} \biggr)
\nonumber
\\
&&\quad\le\Prob\biggl( \biggl\{ \bar{d} \bigl( X_1^n,
\hat{X} \bigl[ \hat{k}_{\mathrm{PML}} \bigl(X_1^n |
\eta\log n\bigr) \bigr]_1^n \bigr) > \frac{\beta_2}{\pinf^{2}}
g_n + \frac{1}{n^{1/2-\mu_n}} \biggr\} \cap G_n \cap
H_n \biggr)
\nonumber
\\
&&\qquad{} + \Prob(\bar{G}_n ) + \Prob(\bar
{H}_n )
\\
&&\quad\le\Prob\biggl( \biggl\{ \bar{d} \bigl( X_1^n,
\hat{X} \bigl[ \hat{k}_{\mathrm{PML}} \bigl(X_1^n |
\eta\log n\bigr) \bigr]_1^n \bigr) > \frac
{\beta_2}{\pinf^{2}}
\bar{\gamma} \bigl( \hat{k}_{\mathrm{PML}} \bigl(X_1^n
| \eta\log n\bigr) \bigr) + \frac{1}{n^{1/2-\mu_n}} \biggr\} \cap H_n
\biggr)
\nonumber
\\
&&\qquad{} + \Prob(\bar{G}_n ) + \Prob(\bar
{H}_n ).\nonumber
\end{eqnarray}
The three terms on the right of (\ref{eqfinal}) is bounded as follows.\vadjust{\goodbreak}

Since the process is non-null with summable continuity rate, Lemma \ref
{lemapprox} in the \hyperref[app]{Appendix} with $\mu=\mu_n$, $\nu\log
n=k_n$ and
$k=\hat{k}_{\mathrm{PML}} (X_1^n | \eta\log n)$ gives
%
%
\begin{eqnarray}
\label{eqf1} &&\mathrm{Pr} \biggl( \biggl\{ \bar{d} \bigl( X_1^n,
\hat{X} \bigl[ \hat{k}_{\mathrm{PML}} \bigl(X_1^n |
\eta\log n\bigr) \bigr]_1^n \bigr) > \frac
{\beta_2}{\pinf^{2}}
\bar{\gamma} \bigl( \hat{k}_{\mathrm{PML}} \bigl(X_1^n
| \eta\log n\bigr) \bigr) + \frac{1}{n^{1/2-\mu_n}} \biggr\} \cap H_n
\biggr)
\nonumber
\\
&&\quad\le2 \RMe ^{1/\RMe } |A|^{k_n+2} \exp\biggl\{ -
\frac{ \pinf^2 }{ 16 \RMe  |A|^3 (\alpha+ \pinf) (\beta_1+1)^2 } \frac{
(n-k_n) 4^{-k_n |\log\pinf|} }{ (1+k_n) n }
\\
&&\hspace*{90pt}{}\times \biggl[ 4^{\mu_n\log n} - \frac{k_n|\log\pinf|(\beta_1+1)^2}{2}
\biggr] \biggr\} .\nonumber
\end{eqnarray}
By Remark~\ref{remgammaalpha}, the summability of the continuity rate
implies the $\alpha$-summability. Hence, for the non-null process with
summable continuity rate and uniformly convergent restricted continuity
rate with parameters $\theta_1$, $\theta_2$, $k_{\theta}$, Proposition
\ref{propqmin} implies that
%
%
\begin{equation}
\label{eqGbound} \Prob(\bar{G}_n ) \le12 \RMe ^{1/\RMe } \exp
\biggl( -\frac{7\alpha_0 (\log|A|)^3 \eta^3}{4\RMe (\alpha+\alpha_0)} \frac
{ n^{\eta2\log|A|} }{ \log n } + \bigl(\eta\log|A|\bigr)\log n \biggr),
\end{equation}
if (\ref{eqnbound}) holds because
\begin{eqnarray*}
&&\Prob\bigl( \bar{\gamma} \bigl( \hat{k}_{\mathrm{PML}} \bigl(X_1^n
| \eta\log n\bigr) \bigr) \ge g_n \bigr)
\\
&&\quad= \Prob\biggl( \hat{k}_{\mathrm{PML}} \bigl(X_1^n
| \eta\log n\bigr) \\
&&\qquad\hspace*{14pt}\le\min\biggl\{ \biggl\lfloor\frac{\eta}{\theta
_2}\log n \biggr
\rfloor,
k\ge0\dvtx\bar{\gamma}(k) < \biggl( \frac{ 6\max( \sqrt{n}
,(|A|-1) \operatorname{pen}(n)) }{ \pinf n^{1-\eta\log
(|A|^4/\pinf)} }
\biggr)^{1/(2\theta_1)} \biggr\} \biggr) .
\end{eqnarray*}
Applying Proposition~\ref{propqminupper} with $r_n=\eta\log n$,
\[
m_n = \min\biggl\{ \lfloor\eta\log n \rfloor, k\ge0\dvtx\bar{
\gamma}(k) < \frac{c \operatorname{pen}(n)}{ n } \biggr\}
\]
and $k_n=h_n+m_n$, it follows that
%
%
\begin{eqnarray}
\label{eqHbound}
&&\Prob(\bar{H}_n ) \le\exp\biggl( -\bigl(|A|-1\bigr)
|A|^{k_n+1} \nonumber\\
&&\hspace*{63pt}{}\times\operatorname{pen}(n) \biggl[ 1 - \frac
{1}{|A|^{1+h_n}} -
\frac{1}{2\operatorname{pen}(n)} \bigl( \log n - k_n \log|A| \bigr) \biggr]
\\
&&\hspace*{63pt}{}+ \frac{c \operatorname{pen}(n)}{\pinf/\log \RMe } + |A|^{k_n+1} C_{\mathrm
{KT}} + \log(\eta\log
n) \biggr) .
\nonumber
\end{eqnarray}
Finally, applying the bounds (\ref{eqf1}), (\ref{eqGbound}) and
(\ref
{eqHbound}) to the right of (\ref{eqfinal}), the proof is
complete.\vspace*{-4pt}
\end{pf*}

\begin{appendix}\label{app}
\section*{Appendix}
%
%
\setcounter{definition}{0}
\begin{lemma}\label{lemsch}
For two probability distributions $P_1$ and $P_2$ on $A^k$,
\[
\bigl\llvert H(P_1) - H(P_2) \bigr\rrvert\le
\frac{1}{\log \RMe } \bigl[ k\log|A| - \log d_{\TV} ( P_1,
P_2 ) \bigr] d_{\TV} ( P_1, P_2 ),
\]
if $d_{\TV} ( P_1, P_2 ) \le1/\RMe $, where
\[
H(P_i) = - \sum_{a_1^k\in A^k} P_i
\bigl(a_1^k\bigr) \log P_i
\bigl(a_1^k\bigr)
\]
is the entropy of $P_i$, $i=1,2$, and
\[
d_{\TV} ( P_1, P_2 ) = \sum
_{a_1^k\in A^k} \bigl\llvert P_1 \bigl(a_1^k
\bigr) - P_2 \bigl(a_1^k\bigr) \bigr\rrvert
\]
is the total variation distance of $P_1$ and $P_2$.
\end{lemma}
\begin{pf}
See Lemma 3.1 of~\cite{Sch}.
\end{pf}
%
%
\begin{lemma}\label{lemktml}
There exists a constant $C_{\mathrm{KT}}$ depending only on $|A|$, such
that for any $0\le k<n$
\[
\log\mathrm{ML}_{k}\bigl(X_1^n\bigr) -\log
P_{\mathrm{KT}, k} \bigl(X_1^n\bigr) \le
C_{\mathrm{KT}} |A|^k + \frac{|A|-1}{2} |A|^k \log
\frac{n}{|A|^k} .
\]
\end{lemma}
\begin{pf}
The bound, see, for example, (27) in~\cite{Cs},
\begin{eqnarray*}
&&\biggl\llvert\log P_{\mathrm{KT}, k} \bigl(X_1^n
\bigr) + k\log|A| -\log\mathrm{ML}_{k}\bigl(X_1^n
\bigr) + \frac{|A|-1}{2} \mathop{\sum_{a_1^k\in
A^k:}}_{N_{n-1}(a_1^k)\ge1}
\log N_{n-1}\bigl(a_1^k\bigr) \biggr\rrvert
\\
&&\quad\le C_{\mathrm{KT}}' |A|^k,
\end{eqnarray*}
where $C_{\mathrm{KT}}'$ depends only on $|A|$, implies the claim using
\[
\mathop{\sum_{a_1^k\in A^k:}}_{N_{n-1}(a_1^k)\ge1} \log
N_{n-1}\bigl(a_1^k\bigr) \le|A|^k
\log\frac{n}{|A|^k} ,
\]
see Proof of Theorem 6 in~\cite{Cs}.
\end{pf}
%
%
\begin{lemma}\label{lemapprox}
Let $X$ be a non-null stationary ergodic process with summable
continuity rate. Then, for any $\mu>0$ and $k\le\nu\log n$, $\nu>0$,
the empirical $k$-order Markov estimator of the process satisfies
\begin{eqnarray*}
&&\mathrm{Pr} \biggl\{ \bar{d} \bigl( X_1^n,
\hat{X}[k]_1^n \bigr) > \beta_2
\pinf^{-2} \bar{\gamma}(k) + \frac{1}{n^{1/2-\mu}} \biggr\}
\\
&&\quad\le2 \RMe ^{1/\RMe } |A|^{2+\nu\log n}
\\
&&\qquad\hspace*{0pt}{}\times \exp\biggl\{ - \frac{ \pinf^2 }{ 16 \RMe  |A|^3 (\alpha+
\pinf) (\beta_1+1)^2 } \frac{
(n-\nu\log n) n^{-2\nu|\log\pinf|} }{ (1+\nu\log n) n }
\\
&&\qquad\hspace*{30pt}{}\times \biggl[ n^{2\mu} - \frac{\nu|\log\pinf|(\beta
_1+1)^2\log n}{2} \biggr] \biggr\}
.
\end{eqnarray*}
\end{lemma}
\begin{pf}
See the proof of Theorem 2 and Lemma 3 in~\cite{CsT3}.
\end{pf}
\end{appendix}

\section*{Acknowledgements}

The author would like to thank the referees for their comments that
helped improving the presentation of the results and generalizing the
consistency concept. The research of the author was supported in part
by NSF Grant DMS-09-06929.


%

\printhistory

\end{document}